\documentclass[12pt,reqno]{amsart}
\usepackage[normalem]{ulem}
\usepackage[a4paper,margin=2.5cm]{geometry}
\usepackage{amsmath,amsthm,amssymb,mathrsfs}
\usepackage{bbm}
\newcommand{\bbm}[1]{\mathbbm{#1}}
\ifdefined\XeTeXversion
  \usepackage{dsfont}
  \renewcommand{\bbm}[1]{\mathds{#1}}
\fi
\usepackage{cite}
\usepackage{xcolor}
\usepackage{etoolbox}
\makeatletter
\patchcmd{\@sect}{\@tocwrite{#1}{#8}}{\@tocwrite{#1}{#7}}{}{\errmessage{First TOC title patch failed}}
\patchcmd{\@sect}{\@tocwrite{#1}{#8}}{\@tocwrite{#1}{#7}}{}{\errmessage{Second TOC title patch failed}}
\makeatother
\usepackage[colorlinks=true]{hyperref}
\hypersetup{
  linkcolor=blue,
  citecolor=green!40!black,
  urlcolor=red
}
\newcommand{\one}{\bbm{1}}

\newtheorem{theorem}{Theorem}[section]

\newtheorem{lemma}[theorem]{Lemma}
\newtheorem{corollary}[theorem]{Corollary}
\newtheorem{example}[theorem]{Example}
\newtheorem*{example*}{Example}
\newtheorem{proposition}[theorem]{Proposition}

\numberwithin{equation}{section}

\newcommand{\Q}{\mathbb Q}
\newcommand{\N}{\mathbb N}
\newcommand{\Z}{\mathbb Z}
\newcommand{\C}{\mathbb C}
\newcommand{\R}{\mathbb R}
\newcommand{\T}{\mathbb T}
\newcommand{\Ebb}{\mathbb E}
\newcommand{\1}{\bbm{1}}

\newcommand{\m}{\mathfrak m}

\newcommand{\Sinf}{\mathbb N_0^{(\mathbb N)}}

\begin{document}
\title[$L^p$-Khinchin sequences]{$L^p$-Khinchin sequences and transference principles II}

\begin{abstract}
A  sequence of distinct positive
integers \(\Lambda=(\lambda_n)_{n\ge1}\) is called an $L^p$-Khinchin sequence  ($1\le p\le \infty$) if for any $f\in L^p(\mathbb{T})$, the Khinchin averages
$\Lambda_nf(x) := \frac1n\sum_{k=1}^n f(\lambda_k x)$ converge to $\int_\mathbb{T} fd \m $ for $\m$-almost all $x$
($\m$ being the Lebesgue measure). 
{The $L^\infty$-Khinchin sequences are characterized by the continuity at $0$ in the topology of convergence in measure of the associated Khinchin maximal operator according to the Bellow-Jones maximal principle.}
We establish an $L^\infty$-transference
principle which allows us to 
compare the Khinchin action with  all  
measure-preserving multiplicative semigroup $(\N, \cdot)$-actions on probability spaces.
As consequences of this principle,  we obtain several practical criteria.  Then we prove 
that {a strictly increasing} $L^\infty$-Khinchin sequence grows at least sub-exponentially, 
and that non-$L^\infty$-Khinchin sequences include: the {non-trivial} translates and the exponentiations
of  sequences in any finitely generated semigroup,  
 almost sure sampling of {nondegenerate} Bernoulli random multiplicative sequences generated by pairwise coprimes, 
 multiplicative sequences arising from some non-balanced substitution sequences etc.
The $L^p$-Khinchin sequences for $1\le p<\infty$  {were} discussed in \cite{FanFanLiLp2026}.
\end{abstract}

\author{Aihua FAN}
\address{LAMFA, UMR 7352 CNRS, University of Picardie, 33 rue Saint Leu, 80039 Amiens, France}
\email{ai-hua.fan@u-picardie.fr}

\author{Shilei FAN}
\address{School of Mathematics and Statistics, and Key Lab NAA--MOE, Central China Normal University, Wuhan 430079, China}
\email{slfan@ccnu.edu.cn}

\author{Tongjin LI}
\address{School of Mathematics and Statistics \& Hubei Key Laboratory of Mathematical Sciences, Central China Normal University, Wuhan, 430079, China}
\email{tongjin.li@mails.ccnu.edu.cn}

\thanks{The authors were supported by NSFC (grants No. 12331004, No.  12231013 and No. 12571205)}

\makeatletter
\@namedef{subjclassname@2020}{\textup{2020} Mathematics Subject Classification}
\makeatother

\maketitle

\tableofcontents

\section{Introduction and main results}
\label{sec:introduction-main-results}

Let \(\Lambda=(\lambda_n)_{n\ge1}\) be a sequence of distinct positive
integers, which is not assumed to be increasing.
Let $\mathbb{T}$ be the additive group $\mathbb{R} / \mathbb{Z}$ 
identified with the unit interval $[0,1)$ and
endowed with Lebesgue measure \(\mathfrak m\).  It is well known that,
for Lebesgue almost every \(x\in\mathbb T\), the sequence
\((\lambda_nx)_{n\ge1}\) is uniformly distributed modulo one.  Hence, 
for every Riemann integrable function
\(f\) on \(\mathbb T\),
\begin{equation}\label{KH-property}
\lim_{n\to\infty} \frac1n\sum_{k=1}^n f(\lambda_k x)
=
\int_{\mathbb T} f\,d\mathfrak m
\qquad
\text{for }\mathfrak m\text{-a.e. }x\in\mathbb T .
\end{equation}
Following
\cite{FFQQ}, we say that \(\Lambda\) is an \emph{\(L^p\)-Khinchin sequence},
for some  \(1\le p\le\infty\), if \eqref{KH-property} holds  for every \(f\in L^p(\mathbb T)\).

In our previous paper \cite{FanFanLiLp2026}, we  carried out a study of
$L^p$-Khinchin sequences for $1\le p<\infty$. 
In 1923, Khinchin \cite{Khintchin} 
conjectured that the full sequence of natural numbers  \((n)_{n\ge1}\) is an
\(L^\infty\)-Khinchin sequence.  This conjecture was disproved  by Marstrand \cite{Marstrand}  in 1970, who 
constructed a counterexample $f \in L^\infty(\mathbb{T})$ which is the indicator function of some  Borel set.
Later,  Bourgain \cite{Bourgain1988} and
Quas--Wierdl \cite{Quas_Wierdl2010} gave different proofs of Marstrand's result.  On
the other hand, for every integer \(q\ge2\), 
 the geometric sequence
\((q^n)_{n\ge1}\) is \(L^1\)-Khinchin by the ergodic theorem (Raikov
\cite{Raikov} and Riesz \cite{Riesz} studied the sequence). 
In \cite{Marstrand}, 
Marstrand also proved that the increasing enumeration of the Furstenberg semigroup
\(\{2^m3^n:m,n\ge0\}\) is an $L^\infty$-Khinchin sequence, and later Nair \cite{Nair1990}
proved the sequence is in fact \(L^1\)-Khinchin using Bewley's ergodic theorem for
$(\mathbb N_0^d,+)$-actions. 

Finding \(L^p\)-Khinchin sequences among subsequences of
\((2^n)_{n\ge1}\) is closely related to the investigation of universal good or bad sequences in ergodic
theory.  Bourgain's results \cite{B1988,Bourgain1989} imply that
\((2^{n^2})_{n\ge1}\) is \(L^p\)-Khinchin for every \(p>1\), whereas  Buczolich and Mauldin's result  \cite{Buczolich2010} implies that it is not
\(L^1\)-Khinchin.  Let us also point out that a result of Rosenblatt \cite{Rosenblatt1991} implies  that
\((2^{2^n})_{n\ge1}\) is not \(L^\infty\)-Khinchin.

In \cite{FanFanLiLp2026}, we characterized $L^p$-Khinchin sequences by a weak $(p,p)$-maximal inequality of the corresponding %
{Khinchin maximal operator}
when $1\le p<\infty$. 
The present paper, a continuation of \cite{FanFanLiLp2026}, is mainly devoted to the 
study of $L^\infty$-Khinchin sequences. As in \cite{FanFanLiLp2026}, we denote
\begin{equation}\label{eq:Lambda_n}
\Lambda_nf(x):=\frac1n\sum_{k=1}^n f(\lambda_kx),
\qquad
\Lambda^*f(x):=\sup_{n\ge1}|\Lambda_n f(x)|
\end{equation}
and call $\Lambda_n$ the {\em Khinchin operators} and $\Lambda^*$ the {\em Khinchin  maximal operator} associated to $\Lambda$. 
It is clear that  Banach's principle and its variants are basic tools for studying  the almost everywhere convergence of $\Lambda_nf(x)$. %
This is the approach we followed in  \cite{FanFanLiLp2026}, using Banach's principle and Sawyer's principle for $1\le p<\infty$.

{\em Notation and Convention.}  Throughout this paper, we adopt the following notation.
For a given probability space $(X, \mathcal{A}, \mu)$, 
we denote by $L^\infty(X)$ the space of essentially bounded measurable functions, 
equipped with the norm  $\|\cdot\|_\infty$, and by $B_\infty(X)$ its closed unit ball.
 The space $L^0(X)$ of measurable functions is always equipped with the topology of   convergence in measure
and {$B_\infty(X)$} is equipped with the induced topology, as a subspace of $L^0(X)$. 
We recall the important fact that the topology of convergence in measure
on \(B_\infty(X)\) or on any bounded set 
in $L^\infty(X)$ is the same as that of $L^1$-convergence (see  \cite{BellowJones1996BanachPrinciple}).
For any operator $L$ which maps $B_\infty(X)$ to a bounded subset of $L^\infty(X)$ 
such that $L(0)=0$, we define its {\em modulus of 
 continuity at $0$}  by
\begin{align}\label{continuity-modulus}
\mathcal D_\delta(L)
:=
\sup\left\{
\|L f\|_1:
 f\in B_\infty(X),\ \|f\|_1<\delta
\right\}, \qquad \delta >0.
\end{align}
Observe that $L$ is continuous at  zero if and only if {$\mathcal D_\delta(L)\to 0$ as $\delta \to 0$.}
We also adopt the notations $$\N= \{ 1,2, 3, \ldots\},\quad  \N_0=\{0,1,2,\ldots \}.$$
{We write $\mathbb N_0^{(\mathbb N)}$ for the direct sum of countably many copies of the additive semigroup $(\mathbb N_0,+)$.}
{Semigroups $\mathbb{N}_0^d$ ($d \ge 1$) are identified as sub-semigroups of $\mathbb{N}_0^{(\mathbb{N})}$.} 
We reserve $\gamma=(\gamma^n)_{n\in\N}$ for multiplicative $(\N,\cdot)$-actions and use $\beta=(\beta^{\mathbf n})$ for additive $(\N_0^d,+)$- or $(\Sinf,+)$-actions. 
We write $\Lambda_n^\gamma$ for averages along $\Lambda$ and $A_n^\beta$ 
for weighted averages of additive actions; the corresponding maximal operators are $\Lambda^{\gamma,*}$ and $A^{\beta,*}$ 
 respectively.
{These notations will be explicitly defined later.}

By applying  the Bellow--Jones continuity
principle (see \cite{BellowJones1996BanachPrinciple})  to the Khinchin operators $(\Lambda_n)_{n \ge 1}$, 
   the following  characterization (\!\!\cite[Theorem 2.2]{FFQQ}) {was proved}: a
sequence of distinct positive integers \(\Lambda\) is \(L^\infty\)-Khinchin if and only if the maximal operator
\[
\Lambda^*:B_\infty(\mathbb T)\longrightarrow L^\infty (\mathbb T)
\]
is continuous at \(0\) (also see Theorem~\ref{prop:Linfty-khinchin-continuity} below). We emphasize that  both spaces $B_\infty(\mathbb{T})$ and $L^\infty(\mathbb{T})$ carry the topology of
convergence in measure. 
In other words,  $\Lambda$ is $L^\infty$-Khinchin if and only if $\lim_{\delta \to 0}\mathcal{D}_\delta(\Lambda^*)=0$.

Our main purpose is to obtain practical criteria for the $L^\infty$-Khinchin property of $\Lambda$, namely  for the continuity at $0$ of $\Lambda^*$. %

The Khinchin problem actually concerns a special action of the multiplicative semigroup of positive integers
$(\mathbb{N}, \cdot)$ on the Lebesgue probability space $(\mathbb{T}, \mathcal{B}(\mathbb{T}), \mathfrak{m})$:
any positive integer $a\in \mathbb{N}$ defines a measure-preserving mapping
\begin{equation}\label{eq:Naction}
     T_a x = ax \mod 1
\end{equation}
and we have $T_a \circ T_b = T_{ab}$. We call this action $(T_a)_{a\in \mathbb{N}}$ the {\em Khinchin action}. 

Our basic idea is to compare  the Khinchin action with other 
{$(\N,\cdot)$-actions} on different probability spaces. 
All these actions are linked by a so-called transference principle (see Theorem~\ref{thm:homomorphism-criterion} below). 
Any such %
{$(\N,\cdot)$-action} will be considered as a model of the Khinchin action 
and it can provide useful information about the Khinchin action.

We even consider general semigroup actions.  Here, by a semigroup we mean a commutative semigroup with at most countably many generators.
 As we shall see,  the Khinchin action   is a representative for all
such semigroup actions.
Let $S$ be a commutative semigroup with identity $e$. A 
measure-preserving $S$-action on a 
{measure} space $(X, \mathcal{A}, \mu)$ 
{which will usually be a probability space},  denoted 
by ${\alpha}: S \curvearrowright (X, \mathcal{A}, \mu)$,  is defined by  a 
family of measure-preserving transformations  ${\alpha}=({\alpha}^s)_{s \in S}$ on $(X, \mathcal{A}, \mu)$ such that 
${\alpha}^{s_1} \circ {\alpha}^{s_2} = {\alpha}^{s_1s_2}$ for all 
$s_1, s_2 \in S$ and ${\alpha}^e$ is the identity for the identity $e \in S$. Here the semigroup operation is multiplicatively denoted.
A measure-preserving  semigroup action ${\alpha}: S \curvearrowright (Y, \mathcal{B}, \sigma) $ 
 is said to be \emph{free} if 
$$
  \sigma\{ y \in Y: {\alpha}^{s_1} y ={\alpha}^{s_2}y\} =0
  $$
  whenever $ \ s_1,  s_2 \in S$  are distinct. It is easy to see that the Khinchin action on $(\mathbb{T}, \mathcal{B}(\mathbb{T}), {\m})$ is free.

For a measure-preserving {$(\N,\cdot)$-action}
${\gamma}=({\gamma^n})_{n \in \mathbb{N}}: {(\N,\cdot)} \curvearrowright (Y, \mathcal{B}, \sigma)$,
let 
\begin{equation}\label{eq:L-beta}
\Lambda_n^{{\gamma}} g(y) := \frac{1}{n} \sum_{k=1}^n g({\gamma}^{\lambda_k} y)
\end{equation}
denote the {\em  averages of ${\gamma}$ along $\Lambda$}. Let $\Lambda^{{\gamma}, *}$ be the corresponding {\em maximal operator} defined by 
\begin{equation}\label{eq:L-beta*}
 \Lambda^{{\gamma}, *} g(y) = \sup_{n\ge 1} |\Lambda_n^{{\gamma}} g(y) |.
\end{equation}
Notice that the Khinchin operator $\Lambda_n$ together with the maximal operator $\Lambda^*$ defined in \eqref{eq:Lambda_n} is a special case
of \eqref{eq:L-beta} and \eqref{eq:L-beta*}. 

Here is our transference characterization of $L^\infty$-Khinchin sequences.

\begin{theorem}\label{thm:homomorphism-criterion}
Let $\Lambda=(\lambda_k)_{k\ge1}$ be a sequence of distinct positive
integers. The following are equivalent.
\begin{enumerate}
\item[\textup{(i)}] $\Lambda$ is an $L^\infty$-Khinchin sequence.
\item[\textup{(ii)}]
$\displaystyle\lim_{{\delta \to 0}}\sup_{{\gamma}}
\mathcal D_\delta(\Lambda^{{\gamma},*})=0$,
where ${\gamma}$ ranges over all measure-preserving multiplicative
{$(\N,\cdot)$-actions} on probability spaces.
\item[\textup{(iii)}]
$\displaystyle\lim_{{\delta \to 0}}
\mathcal D_\delta(\Lambda^{{\gamma},*})=0$ for some free measure-preserving
multiplicative {$(\N,\cdot)$-action} ${\gamma}$ on an atomless standard probability
space.
\end{enumerate}
\end{theorem}

The implication (i) $\Longrightarrow $ (ii)  in Theorem~\ref{thm:homomorphism-criterion} shows that if $\Lambda$ is 
an $L^\infty$-Khinchin sequence, then for any {$(\N,\cdot)$-action} ${\gamma}: {(\N,\cdot)} \curvearrowright (Y, \mathcal{B}, \sigma)$, 
by the Bellow-Jones continuity principle, 
the set of $g \in B_\infty(Y)$
such that $\Lambda_n^{{\gamma}} g(y)$ converges almost everywhere is closed in $B_\infty(Y)$. 
The uniform condition (ii) is often used to establish {the failure of the $L^\infty$-Khinchin property}, as it is easier to violate.
The equivalence (iii) $\Longleftrightarrow $ (i)  in Theorem~\ref{thm:homomorphism-criterion} 
shows that  $\delta \mapsto \mathcal{D}_\delta(\Lambda^{{\gamma}, *})$  is continuous at $0$ for 
some free {$(\N,\cdot)$-action} ${\gamma}$ if and only if $\delta \mapsto \mathcal{D}_\delta(\Lambda^* )$ is continuous at $0$ (namely $\Lambda$ is $L^\infty$-Khinchin).

We shall prove Theorem~\ref{thm:homomorphism-criterion} as a consequence of Theorem
\ref{thm:Conze-Linfty-finite} and Theorem \ref{Conze-N_0^N}, which are transference principles in an abstract form. 
The idea of transference goes back to Calder\'on \cite{Calderon1968}. But a Rokhlin  lemma  for free {$(\mathbb{N}_0^d,+)$-actions} due to Avila and Candela \cite{Avila_Candela2016} is indispensable in our proof of transference principles.

We consider Theorem 
\ref{thm:homomorphism-criterion} 
 as  a Conze-type 
transference principle for {$(\N,\cdot)$-actions} along $\Lambda$.  
Conze \cite{Conze1973} proved  that  the constant involved in the weak $(1,1)$ maximal inequality  
for a free {measure-preserving $\Z$-}action controls the constants for
 all other {measure-preserving $\Z$-}actions.   
Theorem~\ref{thm:homomorphism-criterion} states that 
if the maximal operator for a free {$(\N,\cdot)$-action} is continuous, so is the maximal operator of any action, even  uniformly for all {$(\N,\cdot)$-actions}.
Avila and Candela's Rokhlin  lemma is indispensable for treating non-invertible actions.

Any {$(\N,\cdot)$-action} ${\gamma}: {(\N,\cdot)} \curvearrowright (Y, \mathcal{B}, \sigma)$ 
is considered as a model of the Khinchin action. 
For suitable sequences $\Lambda$ and models ${\gamma}$, it is easier to control the
 behavior of the averaging operators $\Lambda_n^{{\gamma}}$.
This allows us to obtain  information about the Khinchin property of $\Lambda$. 
This is the main advantage  of   Theorem \ref{thm:homomorphism-criterion}.

Based on Theorem \ref{thm:homomorphism-criterion}, we obtain several rather practical criteria  for establishing the 
non-$L^\infty$-Khinchin property.
See \cite{FanFanLiLp2026} for conditions ensuring the $L^p$-Khinchin property.

We first obtain a necessary condition for a strictly increasing sequence to be $L^\infty$-Khinchin:
if a strictly increasing sequence $(\lambda_n)_{n\ge1}$ is
$L^\infty$-Khinchin, then there exists $\eta>0$ such that
\[
 \lambda_n\ge \exp(n^\eta)
 \qquad\text{for all sufficiently large }n.
\]
This can be restated as a growth obstruction to $L^\infty$-Khinchin sequences as follows. 

\begin{theorem}
\label{thm:subexponential-growth-obstruction}
Let $\Lambda=(\lambda_n)_{n\ge1}$ be a strictly increasing sequence
of positive integers. If
\begin{align}\label{subexponential-growth}
   \liminf_{n\to\infty}
 \frac{\log\log\lambda_n}{\log n}=0,
\end{align}
then $\Lambda$ is not an $L^\infty$-Khinchin sequence.
\end{theorem}

Introduce the counting function  
$$\Lambda(x)=\#\{n\ge1:\lambda_n\le x\}.$$
It is easy to see that the condition~\eqref{subexponential-growth} in
Theorem~\ref{thm:subexponential-growth-obstruction} is equivalent to
\begin{equation}\label{subexponential-growth2}
 \forall \ \Delta >0, \qquad \limsup_{x\to\infty}\frac{\Lambda(x)}{(\log x)^\Delta}=\infty.
\end{equation}
Under this condition, Quas and Wierdl~\cite[Corollary~5.6]{Quas_Wierdl2010}
proved that there exists $f\in L^1(\T)$ such that
$\limsup_{N\to\infty}\Lambda_N f(x)=+\infty$ for almost every $x\in\T$. 
{Under the same condition},
 Theorem~\ref{thm:subexponential-growth-obstruction} 
{gives} a bounded function $g$ for which $(\Lambda_N g(x))_{N\ge1}$
diverges on a set of positive measure.

 Theorem~\ref{thm:subexponential-growth-obstruction} will be proved as  a  
 consequence of the following  more general theorem
 which presents an obstruction to $L^\infty$-Khinchin sequences by the existence 
 of a specific semigroup homomorphism into $\mathbb{R}$. 

\begin{theorem}
\label{thm:separated-prefix-concentration}
Let $\Lambda=(\lambda_n)_{n\ge1}$ be a sequence of distinct positive integers. Suppose that there exist a homomorphism
$\ell:(\N,\cdot)\to(\R,+)$ and a constant $\eta>0$ with the following property. For every integer $K\ge1$, there exist $h>0$, real numbers
$t_1,\ldots,t_K$ satisfying $|t_i-t_j|\ge h$ for $i\ne j$, and positive integers $N_1,\ldots,N_K$ such that
\begin{align}\label{cond:concentration}
        \frac1{N_i}\#\{n\le N_i:|\ell(\lambda_n)-t_i|\le h\}
            \ge\eta \qquad(1\le i\le K).
\end{align}
Then $\Lambda$ is not an $L^\infty$-Khinchin sequence.
\end{theorem}

The proof of Theorem~\ref{thm:separated-prefix-concentration} is based 
on Theorem~\ref{thm:homomorphism-criterion}, and is also inspired  by the  
idea of using homomorphisms due to 
Quas--Wierdl \cite{Quas_Wierdl2010}.

In \cite{FFQQ}, the so-called multiplicative sequences were introduced. Given two coprime integers $q_0, q_1 \ge 2$,
any infinite sequence $w \in \{0,1\}^\N$ defines  
$\Lambda(w):=(q_{w_1}\cdots q_{w_n})_{n\ge1}$, called a {\em multiplicative sequence}. 
Question 3 in  \cite{FFQQ}
asks  whether, for every point
$\omega$ of a minimal, uniquely ergodic subshift of $\{2,3\}^{\N}$,
the sequence of partial products $(\omega_1\cdots\omega_n)_{n\ge1}$
is $L^p$-Khinchin.
In \cite{FanFanLiLp2026}, we proved that $\Lambda(w)$ is an $L^1$-Khinchin sequence 
when $w$ is 
$C$-balanced, in particular when $w$ is a
primitive substitutive sequence such that 
the absolute value of the second eigenvalue of the  incidence matrix is strictly less than $1$.
One example is the Fibonacci sequence defined by the substitution $ 0\mapsto 01$ and $1\mapsto 0$. 

We can apply Theorem~\ref{thm:separated-prefix-concentration} to get a negative answer to Question 3 in \cite{FFQQ}. 
Let $w=(w_n)_{n\ge1}$ be the fixed point
of the substitution
\begin{equation}\label{eq:non-balanced-substitution}
{\sigma (0)=010},\qquad\sigma(1)=011.
\end{equation}
Then for coprime integers $q_0, q_1 \ge 2$,
$\Lambda(w)=(q_{w_1}\cdots q_{w_n})_{n\ge1}$ is not
$L^\infty$-Khinchin (Corollary~\ref{counter-example-primitive}).
The incidence matrix of this substitution is
$\left(\begin{smallmatrix}2&1\\1&2\end{smallmatrix}\right)$,
which is strictly positive, so the substitution is primitive and the associated 
subshift is  minimal and uniquely ergodic (see
\cite{Queffelec}). 
Thus, minimality and unique ergodicity of the associated subshift do not imply the $L^\infty$-Khinchin property of the corresponding multiplicative sequence.
However, the Thue--Morse sequence behaves differently: any sequence in the Thue-Morse shift produces an $L^1$-Khinchin sequence
(see \cite{FFQQ}). 

The next homomorphism obstruction to $L^\infty$-Khinchin sequences  is another consequence of Theorem~\ref{thm:homomorphism-criterion}. Its proof is 
inspired by Bourgain's divergence result \cite[Proposition~10]{Bourgain1988}, 
but  uses Karagulyan's generalized result \cite{Karagulyan2011}. In  \cite[Theorem~4.1]{FanFanLiLp2026}, an obstruction involving the behavior of the homomorphism  along the full sequence
was established for $L^p$-Khinchin sequences
for $1\le p<\infty$.  The next theorem provides an obstruction 
to $L^\infty$-Khinchin sequences {based only on} the behavior  along a set of positive upper
density.
 Recall that the {\em  upper density} of a subset of integers $I\subset \mathbb{N}$ is defined by
$$
 \overline{d}(I): = \limsup_{n \to \infty} \frac{\#(I \cap [1,n])}{n}.
 $$

\begin{theorem}\label{Bourgain-obstruction}
  Let $\Lambda=(\lambda_n)_{n \ge 1}$ be a sequence of 
  distinct positive integers. Suppose there exist a homomorphism 
  $\pi: (\N, \cdot) \to (\T, +)$, a point $a \in \T$ and a set $I \subseteq \N$ 
  with positive upper density such that 
  \begin{align}\label{eq:K-cond}
    \lim_{\substack{n \to \infty\\ n \in I}} \pi(\lambda_n) = a, \qquad \pi(\lambda_n) \neq a \ \  \text{for} \ \ n \in I.
  \end{align}
  Then $\Lambda$ is not an $L^\infty$-Khinchin sequence.
\end{theorem}

Theorem~\ref{Bourgain-obstruction} applies when the sequence $(\pi(\lambda_n))$ has only
finitely many subsequential limits and eventually avoids its limit points
(Corollary~\ref{cor-of-Bourgain-obstruction}). Thus any sequence of
distinct positive integers of the form 
\[
{\lambda_n=q^{a_n}r_n \qquad  (q>1, r_n>0, a_n \in \N_0)}
\]
is not $L^\infty$-Khinchin if the sequence $(\log_qr_n\bmod1)$ in $\mathbb{T}$ has only finitely
many subsequential limits and eventually avoids them. Examples include\\
\indent (a)  $(q^n+1)_{n\ge1}$ for integers $q\ge2$ ({Example~\ref{exam:leading-part}});\\
\indent (b)   Every enumeration of
$
\bigcup_{t\in I} \{2^t3^m+1: m\ge0\}
$ for any finite non-empty set $ I\subset\N_0$
 (Example~\ref{prop:union-of-non-K}){.}

Question~6 in \cite{FFQQ} asks whether integer translations preserve the
$L^1$-Khinchin property.  Example~\ref{exam:leading-part}, cited above, 
gives a negative answer. Such a negative answer was
also given by Guo and Zhang~\cite[Corollary~4.4]{GuoZhang2026}.

Theorem \ref{Bourgain-obstruction} also implies that 
for any sequence of distinct primes $(p_n)_{n\ge 1}$, both the sequence itself and its
partial products $(p_1p_2\cdots p_n)_{n\ge 1}$ fail to be $L^\infty$-Khinchin
(Corollary~\ref{cor:prime-subsequence-not-khinchin}). This was also proved in 
\cite[Proposition~4.1 and Corollary~4.2]{GuoZhang2026}.

The following $p$-adic homomorphism obstruction will be deduced from Theorem \ref{Bourgain-obstruction}.

\begin{theorem}
\label{thm:padic-homomorphism-obstruction}
Let $\Lambda=(\lambda_n)_{n\ge1}$ be a sequence
of distinct positive integers. Let $p$ be a prime number.
{Suppose that there exist a semigroup homomorphism $\pi_p: (\N, \cdot) \longrightarrow (\Q_p, +)$, a point 
$a \in \Q_p$ and a subset $I (\subset \N)$ of positive upper density such that }
\begin{equation}\label{eq:padic-concentration}
    \lim_{\substack{n \to \infty\\ n \in I}} \pi_p(\lambda_n) = a, \qquad \pi_p(\lambda_n) \neq a \ \text{for every} \ n \in I.
\end{equation}
Then $\Lambda$ is not an $L^\infty$-Khinchin sequence.
\end{theorem}

Rosenblatt's theorem \cite[Theorem 23]{Rosenblatt1991} implies that $(2^{2^n})$ is not $L^\infty$-Khinchin.
Let $v_2$ denote the $2$-adic valuation.
{Using} Theorem \ref{thm:padic-homomorphism-obstruction} with the homomorphism $m\mapsto v_2(m)$ provides a new proof
 of Rosenblatt's result (see Corollary \ref{cor:double-exponential-padic}).
 
Let $S(q_1,\ldots,q_d)$ be the semigroup generated by pairwise coprime integers
$q_1,\ldots q_d \ge 2$ and {let} $P$ be a nonconstant polynomial with $P(\N) \subseteq\N$ and $P(0)\ne0$.
Then every sequence of distinct elements of $P(S(q_1, \ldots, q_d))$ is not $L^\infty$-Khinchin. 
In particular, any sequence of distinct integers in $S(q_1, \ldots, q_d)+1$ is 
not $L^\infty$-Khinchin.

Fix  a sequence of pairwise coprime
positive integers $(q_j)_{j \ge 1}$ with $q_j \ge 2$.
Each infinite sequence $\omega=(\omega_n)_{n \ge 1} \in \N^{\N}$ produces a multiplicative sequence 
$\Lambda(\omega)=(\lambda_n)$ defined by 
\begin{align*}
  \lambda_n=q_{\omega_1}q_{\omega_2} \cdots q_{\omega_n}, \qquad n \ge 1.
\end{align*}
We will select $\omega_j$ independently according to a 
probability measure $\rho$ on $\mathbb{N}$.

\begin{theorem}\label{thm:random-sequence}
Let $\rho$ be a probability measure on $\N$ with $\# \operatorname{supp}{\rho} \ge 2$. 
Let $\rho^\mathbb{N}$ denote the associated Bernoulli product probability on $\N^{\N}$.  
Then for $\rho^{\N}$-almost all $\omega \in \N^\N$, the random sequence $\Lambda(\omega)$ is 
not an $L^\infty$-Khinchin sequence.
\end{theorem}

The proof of Theorem \ref{thm:random-sequence} is based on  Theorem \ref{thm:homomorphism-criterion} (ii) and on a 
construction inspired by
Lacey--Petersen--Wierdl--Rudolph (\!\cite{Lacey1994}, Theorem~4 and the remarks following it) using  Strassen's invariance principle.

We shall first present 
the maximal principles in Section 
~\ref{sec:transference-preliminaries} and prove the transference principle 
in Section~\ref{sec:Linfty-transference}. Then we shall deduce the above-mentioned results almost following the order of their appearance 
in this introduction. 


\section{Characterization of $L^\infty$-Khinchin sequences by the Bellow--Jones principle}
\label{sec:transference-preliminaries}

 We first recall the Bellow--Jones criterion, an $L^\infty$-counterpart of Sawyer's principle, which concerns weak-type $(p,p)$ maximal inequalities in the $L^p$ setting
 when $1\le p<\infty$. %
 {We then apply it to Khinchin averages.}

\subsection{$L^\infty$-Banach principle due to Bellow--Jones}

Consider a probability space $(X, \mathcal{B}, \mu)$. Recall that 
\(L^0(X)\) denotes  the space of measurable functions and we always  equip \(L^0(X)\) with the topology of
convergence in measure, which is  metrizable with the metric 
\[
d_\mu(f,g):=\int_X\min\{|f-g|,1\}\,d\mu.
\]
We  denote by $B_\infty(X)$ the closed unit ball of the space \(L^\infty(X)\) of {essentially} bounded {measurable} functions:
\[
B_\infty(X):=\{f\in L^\infty(X):\|f\|_\infty\le1\}.
\]
Both  \(L^\infty(X)\) and $B_\infty(X)$ are always considered as topological subspaces of 
\(L^0(X)\).  
On $B_\infty(X)$, and more generally on every bounded subset of $L^\infty(X)$, 
the topology of convergence in measure coincides with the topology induced by 
any $L^p$ norm, $1\le p<\infty$
 (see \cite[Proposition~1]{BellowJones1996BanachPrinciple}).

We use the following  Bellow--Jones Banach
principle.  

\begin{theorem}[Bellow--Jones, {\cite[Theorem 1, Theorem 2]{BellowJones1996BanachPrinciple}}] 
\label{thm:Bellow-Jones}
Let \((X,\mathcal B,\mu)\) be a probability space, and let
\(L_n:L^\infty(X)\to L^0(X)\), \(n\ge1\), be continuous linear
operators.  Let  \(L^*f:=\sup_{n\ge1}|L_nf|\) be the maximal function of $f\in L^\infty(X)$.
\begin{enumerate}
\item[\textup{(i)}]  If \(L_nf(x)\) converges \text{a.e.} for every
\(f\in L^\infty(X)\), then \(L^*f<\infty\) \text{a.e.}  for every
 \(f\in L^\infty\) and \(L^*:B_\infty(X)\to L^0(X)\) is continuous at \(0\).

\item[\textup{(ii)}] If \(L^*f<\infty\) \text{a.e.}  for every
\(f\in L^\infty(X)\) and if
\(L^*:B_\infty(X)\to L^0(X)\) is continuous at \(0\), then
\[
\bigl\{f\in B_\infty(X):L_nf(x)\text{ converges for almost every }x\bigr\}
\]
is closed in \(B_\infty(X)\). 
\end{enumerate}
\end{theorem}

The above assertion (i) says that the 
almost-everywhere convergence forces the continuity at $0$ of the maximal
operator. 
The assertion (ii) is the closedness mechanism. 

\subsection{Characterization of $L^\infty$-Khinchin sequences}

The following characterization
and its proof are recalled from \cite[Theorem~2.2]{FFQQ}.

\begin{theorem}
  \label{prop:Linfty-khinchin-continuity}
Let \(\Lambda=(\lambda_n)\) be a sequence of distinct positive integers.
Then \(\Lambda\) is an \(L^\infty\)-Khinchin sequence if and only if the
maximal Khinchin operator
\[
\Lambda^*:B_\infty(\mathbb T)\longrightarrow  L^0(\mathbb T)
\]
is continuous at \(0\).
\end{theorem}

\begin{proof}
For completeness, we include here the proof taken from \mbox{\cite{FFQQ}}.

Assume  that \(\Lambda\) is an \(L^\infty\)-Khinchin sequence. 
Then {each $\Lambda_N$ is continuous for convergence in measure and} 
\(\Lambda_nf(x)\) converges a.e. for every \(f\in L^\infty(\mathbb T)\).
By Theorem~\ref{thm:Bellow-Jones} \textup{(i)}, \(\Lambda^*:B_\infty(\mathbb T)\to
L^0(\mathbb T)\) is continuous at \(0\).

Conversely, assume that $\Lambda^*:B_\infty(\mathbb T)\to L^0(\mathbb T)$ is continuous at $0$.
Clearly
\(\|\Lambda^* f\|_{L^\infty(\mathbb T)}
\le \|f\|_{L^\infty(\mathbb T)}\) for all  \(f\in L^\infty(\mathbb T)\). Hence Theorem~\ref{thm:Bellow-Jones}
\textup{(ii)} implies that   the set of functions $f$ in
\(B_\infty(\mathbb T)\) for which \(\Lambda_Nf\) converges \text{a.e.} 
 is closed.  However, since $(\lambda_n x)$ is uniformly distributed for almost all $x$, 
 \(\Lambda_Nf\) converges \text{a.e.}  for every continuous function $f$. 
 Since the continuous functions in $B_\infty(\mathbb T)$ are dense in
$B_\infty(\mathbb T)$ for the topology of convergence in measure, and
the convergence set is closed, it follows that 
\(\Lambda_Nf\) converges \text{a.e.}  for every $f\in B_\infty(\mathbb{T})$, hence for every $f\in L^\infty(\mathbb{T})$.

It remains to identify the limit. Since \(f\in L^\infty(\mathbb T)\subset
L^1(\mathbb T)\), the \(L^1\)-mean convergence theorem
\cite[Theorem~3.1]{FFQQ} gives
\[
\Lambda_nf\longrightarrow\int_{\mathbb T}f\,d\mathfrak m
\qquad\text{in }L^1(\mathbb T).
\]
Thus the almost-everywhere  limit of $\Lambda_nf$ must be the  constant
$\int_{\mathbb T}f\,d\mathfrak m$. Hence \(\Lambda\) is an
\(L^\infty\)-Khinchin sequence.
\end{proof}

Theorem~\ref{prop:Linfty-khinchin-continuity} reduces the Khinchin problem to continuity at $0$ of the maximal operator $\Lambda^*$.
In the sequel, we shall seek practical criteria for this  continuity of $\Lambda^*$.

Let us point out that the Khinchin action is a measure-preserving {$(\Sinf,+)$-action}, where $\Sinf$ is the infinite direct sum of the additive semigroup
$\mathbb{N}_0$. Indeed, 
let
\(\mathfrak p_1<\mathfrak p_2<\cdots\) be the primes. By the fundamental theorem of arithmetic, the multiplicative semigroup $\mathbb{N}$ is isomorphic to the additive semigroup 
$\Sinf$. In fact,  we have the
semigroup isomorphism $\Phi:(\mathbb N,\cdot)\longrightarrow{(\Sinf,+)}$ defined by
\begin{equation}\label{eq:Phi}
\Phi(n):=(v_j(n))_{j\ge1} \quad \text{for} \quad n=\prod_j\mathfrak p_j^{v_j(n)}.
\end{equation}
We call $\Phi$ the {\em lifting mapping}.  This isomorphism $\Phi$ shows that the Khinchin action is nothing but an {$(\Sinf,+)$-action}.

In the next section, we 
will establish an $L^\infty$-transference principle for {$(\Sinf,+)$-actions} and prove 
Theorem \ref{thm:homomorphism-criterion}.

\section{$L^\infty$-transference principle}
\label{sec:Linfty-transference}

The purpose of this section is to prove 
an $L^\infty$-transference principle for {$(\N_0^d,+)$-actions} 
(Theorem \ref{thm:Conze-Linfty-finite}) and 
for   {$(\Sinf,+)$-actions}  (Theorem~\ref{Conze-N_0^N}). Then we 
prove
Theorem \ref{thm:homomorphism-criterion} as a consequence.

Here is the setting in which we work.
Fix a measure-preserving {$(S,+)$-action} $\alpha : {(S,+)} \curvearrowright (X, \mathcal{A}, \mu)$, 
{where $S$ is $\N_0^d$ or $\Sinf$, equipped with componentwise addition.}
Given a sequence of finitely supported weights (i.e., finitely supported functions) $\nu= (\nu_j)_{j \ge 1}$
on $S$, 
we define the corresponding weighted averaging operators $A_j^\alpha$ by
\begin{equation*}
A_j^\alpha f(x)
:=\sum_{\mathbf n\in  S}
\nu_j(\mathbf n)f(\alpha^{\mathbf n}x) \qquad (f: X \to \C).
\end{equation*}

The dependence on the weight sequence $\nu$ is suppressed from the notation, since $\nu$ will remain fixed throughout.
Thus, the behavior of the averaging operators is encoded by the sequence of weights.
The associated  maximal operators are defined by
\[
A^{\alpha,*} f(x)= \sup_{j \ge 1}|A_j^\alpha f(x)|, \qquad 
A_J^{\alpha,*}f(x):=\sup_{1\le j\le J}|A_j^\alpha f(x)| \quad ( J \ge 1).
\]
We always assume that 
\begin{align}\label{boundedness-of-weights}
  0<C_\nu:=\sup_{j \ge 1}\sum_{\mathbf{n} \in S} |\nu_j(\mathbf{n})| < \infty.
\end{align}

We shall study the continuity in measure at $0$ of the maximal operator $A^{\alpha,*}$.
Under the assumption \eqref{boundedness-of-weights}, $A^{\alpha,*}$ maps 
$B_\infty(X)$ to a bounded subset of $L^\infty(X)$. 
Thus, using the notation \eqref{continuity-modulus}, the
continuity at $0$ of $A^{\alpha,*}: B_\infty(X) \to L^0(X)$ 
is characterized by $\lim_{\delta \to 0}\mathcal{D}_\delta(A^{\alpha,*})=0$.

We shall use the following Rokhlin lemma due to Avila--Candela \cite{Avila_Candela2016} which is valid for commuting, possibly noninvertible,
measure-preserving transformations.

\begin{lemma}[Avila--Candela, {\cite[Theorem~1.2]{Avila_Candela2016}}]
\label{lem:Avila-Candela-tower}
Let $\alpha:{(\N_0^d,+)}\curvearrowright(X,\mathcal A,\mu)$ be a free
measure-preserving action on an atomless standard probability space.
For every integer $L\ge1$ and every $\varepsilon>0$, there exists a
measurable set $E\subset X$ such that the sets
\begin{align*}
(\alpha^{\mathbf q})^{-1}E,
\qquad \mathbf q\in Q_L:=\{0,1,\ldots,L-1\}^d,
\end{align*}
are pairwise disjoint and their union has measure greater than
$1-\varepsilon$.
\end{lemma}

\subsection{\texorpdfstring{$L^\infty$}{L-infinity}-transference principle for \texorpdfstring{{$(\mathbb N_0^d,+)$-actions}}{finite-dimensional additive actions}}

\begin{theorem}[$L^\infty$-transference]
\label{thm:Conze-Linfty-finite}
Let $\alpha$ be a free measure-preserving {$(\N_0^d,+)$-action} on an atomless
standard probability space $(X,\mathcal A,\mu)$, and let $(\nu_j)$ be a sequence of weights 
satisfying \eqref{boundedness-of-weights}. The following are equivalent.
\begin{enumerate}
\item[\textup{(i)}] $A^{\alpha,*}:B_\infty(X)\to L^0(X)$ is continuous
at $0$.
\item[\textup{(ii)}]
$\displaystyle\lim_{{\delta \to 0}}\sup_\beta
\mathcal D_\delta(A^{\beta,*})=0$,
where $\beta$ ranges over all measure-preserving {$(\N_0^d,+)$-actions} on
probability spaces.
\item[\textup{(iii)}]
$\displaystyle\lim_{{\delta \to 0}}\sup_{M\ge1}
\mathcal D_\delta(A^{\tau_M,*})=0$,
where $\tau_M^{\mathbf n}y=y+\mathbf n\pmod M$ on
$\Z_M^d:=(\Z/M\Z)^d$, equipped with normalized counting measure.
\end{enumerate}
Moreover, for every such action $\beta$, every $J\ge1$, and every
$\delta,\eta>0$,
\begin{equation}\label{est:contin-modulus}
\mathcal D_\delta(A_J^{\beta,*})
\le 2\mathcal D_\eta(A_J^{\alpha,*})
+\frac{2C_\nu\delta}{\eta}.
\end{equation}
The constants in this estimate are independent of $J$, $d$, and $\beta$.
\end{theorem}

\begin{proof}
Both \textup{(ii)}$\Rightarrow$\textup{(i)} and
\textup{(ii)}$\Rightarrow$\textup{(iii)} are immediate.
We prove the other two implications.

\emph{Proof of \textup{(i)}$\Rightarrow$\textup{(ii)} and of
\eqref{est:contin-modulus}.}
Fix $J\ge1$ and set
\begin{align*}
K_J&:=\bigcup_{j=1}^J\operatorname{supp}\nu_j,
& R_J&:=\max\bigl(\{0\}\cup
\{\|\mathbf n\|_\infty:\mathbf n\in K_J\}\bigr).
\end{align*}
Choose $M$ large enough that
\begin{equation}\label{eq:R+K}
\frac34\frac{M^d}{(M+R_J)^d}>\frac12.
\end{equation}
Set $L=M+R_J$ and $\mathbf p=(L-1,\ldots,L-1)$.
Lemma~\ref{lem:Avila-Candela-tower}, with tower error $1/4$, gives a
base $E$ such that
\begin{align*}
B_{\mathbf r}:=(\alpha^{\mathbf p-\mathbf r})^{-1}E,
\qquad \mathbf r\in Q_L,
\end{align*}
are pairwise disjoint and $\mu(E)|Q_L|>3/4$. Hence
\begin{equation}\label{eq:Conze-interior-mass}
\mu(E)|Q_M|>\frac12.
\end{equation}

Let $\beta:{(\N_0^d,+)}\curvearrowright(Y,\mathcal B,\sigma)$ be
measure-preserving, and take $g\in B_\infty(Y)$ with $\|g\|_1<\delta$.
Choose a representative of $g$ bounded in absolute value by $1$
 and, for $y\in Y$, put
\begin{align*}
h_y(x):=\sum_{\mathbf r\in Q_L}
g(\beta^{\mathbf r}y)\one_{B_{\mathbf r}}(x).
\end{align*}
The disjointness of the levels gives $h_y\in B_\infty(X)$.
If $\mathbf m\in Q_M$, $x\in B_{\mathbf m}$, and $\mathbf n\in K_J$,
then $\mathbf m+\mathbf n\in Q_L$ and
\begin{align*}
\alpha^{\mathbf p-(\mathbf m+\mathbf n)}(\alpha^{\mathbf n}x)
=\alpha^{\mathbf p-\mathbf m}x\in E.
\end{align*}
Thus $\alpha^{\mathbf n}x\in B_{\mathbf m+\mathbf n}$, which yields
\begin{align*}
A_j^\alpha h_y(x)
&=\sum_{\mathbf n}\nu_j(\mathbf n)g(\beta^{\mathbf m+\mathbf n}y)
=A_j^\beta g(\beta^{\mathbf m}y),\qquad 1\le j\le J,\\
A_J^{\alpha,*}h_y(x)&=A_J^{\beta,*}g(\beta^{\mathbf m}y).
\end{align*}
Integrating over $y$, then over the levels $B_{\mathbf m}$ with
$\mathbf m\in Q_M$, gives
\begin{equation}\label{eq:Conze-comparison}
\mu(E)|Q_M|\|A_J^{\beta,*}g\|_1
\le\int_Y\|A_J^{\alpha,*}h_y\|_1\,d\sigma(y).
\end{equation}
Here we used Fubini's theorem and preservation of $\sigma$ by $\beta$.
Similarly,
\begin{align*}
\int_Y\|h_y\|_1\,d\sigma(y)
=\mu(E)|Q_L|\|g\|_1\le\|g\|_1.
\end{align*}
For $\eta>0$, let $G=\{y:\|h_y\|_1<\eta\}$.
On $G$, the integrand in \eqref{eq:Conze-comparison} is at most
$\mathcal D_\eta(A_J^{\alpha,*})$; on $Y\setminus G$, it is at most
$C_\nu$. Markov's inequality therefore gives
\begin{align*}
\int_Y\|A_J^{\alpha,*}h_y\|_1\,d\sigma(y)
\le\mathcal D_\eta(A_J^{\alpha,*})+
\frac{C_\nu\|g\|_1}{\eta}.
\end{align*}
Combining this with \eqref{eq:Conze-interior-mass} and
\eqref{eq:Conze-comparison}, and taking the supremum over $g$, proves
\eqref{est:contin-modulus}.

Under \textup{(i)}, given $\varepsilon>0$, choose $\delta_0>0$ with
$\mathcal D_{\delta_0}(A^{\alpha,*})<\varepsilon/4$, and then take
$0<\delta<\varepsilon\delta_0/(4C_\nu)$.
Since $A_J^{\alpha,*}\le A^{\alpha,*}$, \eqref{est:contin-modulus}
gives a bound strictly smaller than $\varepsilon$, uniformly in
$J$ and $\beta$. Finally, monotone convergence and interchange of
the two suprema give
\begin{align*}
\mathcal D_\delta(A^{\beta,*})
=\sup_{J\ge1}\mathcal D_\delta(A_J^{\beta,*}),
\end{align*}
which proves \textup{(ii)}.

\emph{Proof of \textup{(iii)}$\Rightarrow$\textup{(ii)}.}
Fix $\varepsilon>0$. By \textup{(iii)}, there exists $\delta_0>0$ such that
\begin{equation}\label{eq:Ec-M}
\sup_{M\ge1}\mathcal D_{\delta_0}(A^{\tau_M,*})<\frac\varepsilon3.
\end{equation}
Set
\begin{equation}\label{eq:dd}
\delta:=\frac{\varepsilon\delta_0}{3C_\nu}.
\end{equation}
Take any measure-preserving action
$\beta:{(\N_0^d,+)}\curvearrowright(Y,\mathcal B,\sigma)$ and
$f\in B_\infty(Y)$ with $\|f\|_1<\delta$.
Fix $J$ and use $K_J,R_J$ as above. For $M>R_J$ and $y\in Y$, define
a function on $\Z_M^d$, identified with $Q_M$, by
\begin{align*}
f_y(\mathbf z):=f(\beta^{\mathbf z}y),\qquad\mathbf z\in Q_M.
\end{align*}
We may take $|f|\le1$ everywhere, so $\|f_y\|_\infty\le1$.
If $\mathbf z\in Q_{M-R_J}$, then $\mathbf z+K_J\subset Q_M$.
No reduction modulo $M$ occurs in the first $J$ averages, and hence
\begin{align*}
A_j^{\tau_M}f_y(\mathbf z)
=\sum_{\mathbf n}\nu_j(\mathbf n)f(\beta^{\mathbf z+\mathbf n}y)
=A_j^\beta f(\beta^{\mathbf z}y),\qquad 1\le j\le J.
\end{align*}
In particular,
\begin{equation}\label{transferrence-compare}
A_J^{\beta,*}f(\beta^{\mathbf z}y)
\le A^{\tau_M,*}f_y(\mathbf z),\qquad\mathbf z\in Q_{M-R_J}.
\end{equation}
Preservation of $\sigma$ gives
\begin{equation}\label{eq:HM-average}
\int_Y\|f_y\|_{L^1(\Z_M^d)}\,d\sigma(y)
=\int_Y\frac1{M^d}\sum_{\mathbf z\in Q_M}
|f(\beta^{\mathbf z}y)|\,d\sigma(y)=\|f\|_1.
\end{equation}
Let $G_M=\{y:\|f_y\|_{L^1(\Z_M^d)}<\delta_0\}$.
By Markov's inequality,
\begin{equation}\label{eq:bad-x}
\sigma(Y\setminus G_M)\le\frac{\|f\|_1}{\delta_0}
<\frac{\varepsilon}{3C_\nu}.
\end{equation}
On $G_M$, \eqref{eq:Ec-M} bounds $\|A^{\tau_M,*}f_y\|_1$ by
$\varepsilon/3$. On its complement we use the uniform bound
$\|A^{\tau_M,*}f_y\|_1\le C_\nu$. Consequently,
\begin{equation}\label{eq:integrated-discrete-maximal}
\int_Y\|A^{\tau_M,*}f_y\|_1\,d\sigma(y)
\le\frac\varepsilon3+C_\nu\sigma(Y\setminus G_M)
<\frac{2\varepsilon}{3}.
\end{equation}
Sum \eqref{transferrence-compare} over $\mathbf z\in Q_{M-R_J}$ and
integrate in $y$. Using invariance once more, we obtain
\begin{align*}
\frac{|Q_{M-R_J}|}{M^d}\|A_J^{\beta,*}f\|_1
&=\int_Y\frac1{M^d}\sum_{\mathbf z\in Q_{M-R_J}}
A_J^{\beta,*}f(\beta^{\mathbf z}y)\,d\sigma(y)\\
&\le\int_Y\|A^{\tau_M,*}f_y\|_1\,d\sigma(y)
<\frac{2\varepsilon}{3}.
\end{align*}
Letting $M\to\infty$ and then $J\to\infty$ gives
$\|A^{\beta,*}f\|_1\le2\varepsilon/3$.
This bound is uniform in $\beta$ and in $f$ with $\|f\|_1<\delta$,
so \textup{(ii)} follows.
\end{proof}

{\bf Remark.}
Theorem \ref{thm:Conze-Linfty-finite}  makes continuity at $0$ independent of the chosen free
action and identifies it with a uniform condition on all models.
The proof of the  implication \textup{(i)}$\Rightarrow$\textup{(ii)} uses the
Avila--Candela tower theorem, while the proof of the implication
\textup{(iii)}$\Rightarrow$\textup{(ii)} is a Calder\'on transference
argument. We use {the term ``Conze transference principle''} to indicate the link to Conze's result
for {$(\Z,+)$-actions} \cite{Conze1973}, which motivated the investigation of the transference principle
stated in Theorem \ref{thm:Conze-Linfty-finite} and in Theorem~\ref{Conze-N_0^N} below.

\subsection{\texorpdfstring{$L^\infty$}{L-infinity}-transference principle for \texorpdfstring{{$(\Sinf,+)$-actions}}{infinite-dimensional additive actions}}
We further extend Theorem \ref{thm:Conze-Linfty-finite} to
the case of $\Sinf$. %

\begin{theorem}\label{Conze-N_0^N} 
Let $\alpha$ be a free measure-preserving {$(\Sinf,+)$-action} 
  on an atomless standard probability space $(X, \mathcal{A}, \mu)$.
  Let $(\nu_j)_{j \ge 1}$ be a sequence of finitely
  supported weights on $\Sinf$ satisfying \eqref{boundedness-of-weights}. The following are equivalent:
  \begin{enumerate}
\item[\textup{(i)}] $A^{\alpha, *} : B_\infty(X) \to L^0(X)$ is continuous at $0$.
\item[\textup{(ii)}]   $\lim_{{\delta \to 0}} \sup_{\beta} \mathcal{D}_\delta(A^{\beta,*})=0$,
where $\beta$ ranges over all measure-preserving $(\Sinf,+)$-actions on 
probability spaces.
  \end{enumerate}
\end{theorem}

\begin{proof}
The implication \textup{(ii)}$\Rightarrow$\textup{(i)} is immediate, since
$\alpha$ is among the actions considered in \textup{(ii)}. We prove the
converse.

Fix $J\ge1$, and let
\[
K_J:=\bigcup_{j=1}^{J}\operatorname{supp}\nu_j,
\qquad
F_J:=\{1\}\cup\bigcup_{\mathbf n\in K_J}
\operatorname{supp}\mathbf n.
\]
Thus $F_J$ is finite and nonempty. Write
$d=\#F_J$ and identify $\mathbb N_0^{F_J}$ with $\mathbb N_0^d$.
Extension by zero outside $F_J$ gives an injection
$\iota_J:\mathbb N_0^{F_J}\to\Sinf$.

Given any measure-preserving {$(\Sinf,+)$-action} $\beta$ on  {a probability space}
$(Y,\mathcal B,\sigma)$, define {$(\mathbb N_0^{F_J},+)$-actions} by
\[
\alpha_J^{\mathbf u}:=\alpha^{\iota_J(\mathbf u)},
\qquad
\beta_J^{\mathbf u}:=\beta^{\iota_J(\mathbf u)}.
\]
The action $\alpha_J$ is free because $\alpha$ is free. For $1\le j\le J$, put
$\nu_{j,J}(\mathbf u):=\nu_j(\iota_J(\mathbf u))$.
The first $J$ averages associated with
$(\alpha_J,\nu_{j,J})$ and $(\beta_J,\nu_{j,J})$
are exactly the first $J$ averages associated with
$(\alpha,\nu_j)$ and $(\beta,\nu_j)$, respectively.

Apply the quantitative estimate \eqref{est:contin-modulus} from
Theorem~\ref{thm:Conze-Linfty-finite} to these first $J$ weights,
extended by zero for $j>J$. If all first $J$ weights vanish, the
following bound is immediate. Otherwise, their $\ell^1$ norms are
bounded by $C_\nu$, so for every $\delta_0,\delta>0$,
\begin{equation}\label{eq:infinite-transference-bound}
\mathcal D_\delta(A_J^{\beta,*})
\le
2\mathcal D_{\delta_0}(A_J^{\alpha,*})
+\frac{2C_\nu\delta}{\delta_0}.
\end{equation}
Here the constants are independent of $J$, $d$, and $\beta$.

Since $A_J^{\alpha,*}\le A^{\alpha,*}$, the bound
\eqref{eq:infinite-transference-bound} implies
\[
\sup_\beta\mathcal D_\delta(A_J^{\beta,*})
\le
2\mathcal D_{\delta_0}(A^{\alpha,*})
+\frac{2C_\nu\delta}{\delta_0}
\qquad(J\ge1).
\]
By monotone convergence,
\[
\mathcal D_\delta(A^{\beta,*})
=
\sup_{J\ge1}\mathcal D_\delta(A_J^{\beta,*}).
\]
Hence
\[
\sup_\beta\mathcal D_\delta(A^{\beta,*})
\le
2\mathcal D_{\delta_0}(A^{\alpha,*})
+\frac{2C_\nu\delta}{\delta_0}.
\]
Under \textup{(i)}, the continuity at $0$ of $A^{\alpha,*}:
B_\infty(X) \to L^0(X)$ shows that
$\mathcal D_{\delta_0}(A^{\alpha,*})\to0$ as ${\delta_0 \to 0}$.
First letting ${\delta \to 0}$ with $\delta_0$ fixed, and then
letting ${\delta_0 \to 0}$, proves \textup{(ii)}.
\end{proof}

\subsection[Proof of Theorem \ref{thm:homomorphism-criterion}: Criteria for \texorpdfstring{$L^\infty$}{L-infinity}-Khinchin sequences]{Proof of Theorem \ref{thm:homomorphism-criterion}: Criteria {for} $L^\infty$-Khinchin sequences}

\begin{proof}[Proof of Theorem~\ref{thm:homomorphism-criterion}]
The proof amounts to identifying multiplicative {$(\N,\cdot)$-actions} with {$(\mathbb N_0^{(\mathbb N)},+)$-actions}, and proceeds in several steps.

{\em Step 1. The Khinchin action $(T_a)_{a\in \mathbb{N}}$ is free.}  Indeed, for two distinct positive integers $a$ and $b$, the set
$$
\{x\in \mathbb{T}: T_ax= T_b x\} = \{x\in\T:(a-b)x=0\pmod1\} = \left\{ \frac{j}{|b-a|}: 0 \le j <|b-a|\right\}
$$ is finite, hence of zero Lebesgue measure.

{\em Step 2. Correspondence.} Let $\mathcal{F}_\mathbb{N}$ be the family of all {$(\N,\cdot)$-actions} 
on probability spaces
and let $\mathcal{F}$ be the family of all {$(\mathbb{N}_0^{(\mathbb{N})},+)$-actions} on probability spaces. 
Let $\Phi: {(\mathbb{N},\cdot)} \to {(\mathbb{N}_0^{(\mathbb{N})},+)}$ be the 
 isomorphism  defined by \eqref{eq:Phi}. 
To the given sequence $\Lambda=(\lambda_k)_{k\ge 1}$, we associate the sequence of weights $\nu=(\nu_N)_{N\ge 1}$ where 
$$
\nu_N(\cdot) :=\frac1N \sum_{k=1}^N \delta_{\Phi(\lambda_k)}(\cdot).
$$
{Each multiplicative action $\gamma=(\gamma^n)_{n\in\N}$ in $\mathcal F_\mathbb N$ corresponds uniquely to an additive action $\beta=(\beta^{\mathbf n})_{\mathbf n\in\Sinf}$ in $\mathcal F$, defined by
$\beta^{\Phi(n)}=\gamma^n$. For the above weights $\{\nu_N\}$, the corresponding averages satisfy $\Lambda_n^\gamma=A_n^\beta$. Hence their maximal operators and moduli of continuity agree:}
\begin{equation}\label{eq:Correspondence}
{\Lambda^{\gamma,*}=A^{\beta,*},\qquad
\mathcal D_\delta(\Lambda^{\gamma,*})=\mathcal D_\delta(A^{\beta,*}).}
\end{equation}
{This correspondence preserves freeness, since $\Phi$ is an isomorphism.}

{\em Step 3. Proof of (iii) $\Longrightarrow (ii)$.} Assume $\displaystyle\lim_{{\delta \to 0}}
\mathcal D_\delta(\Lambda^{{\gamma},*})=0$ for some free measure-preserving
multiplicative {$(\N,\cdot)$-action} ${\gamma}$. %
{Let $\beta$ be its { corresponding $\mathbb{N}_0^{(\mathbb{N})}$-action.}} By the correspondence \eqref{eq:Correspondence},
$$\displaystyle\lim_{{\delta \to 0}}
\mathcal D_\delta(A^{\beta,*})=0.
$$
Since the corresponding action remains free, Theorem~\ref{Conze-N_0^N} yields
$$
\lim_{{\delta \to 0}} \sup_{{\beta} \in \mathcal{F}} \mathcal{D}_\delta(A^{{\beta},*})=0.
$$
Again, by the  correspondence \eqref{eq:Correspondence} but used in the other direction, we get 
$$
\lim_{{\delta \to 0}} \sup_{{\gamma} \in \mathcal{F}_\mathbb{N}} \mathcal{D}_\delta(\Lambda^{{\gamma},*})=0.
$$
This is the statement of (ii).

{\em Step 4. Proof of equivalences.} 
By Theorem~\ref{prop:Linfty-khinchin-continuity}, (i) means that the Khinchin maximal operator $\Lambda^*$ is continuous at $0$, namely 
$\lim_{{\delta \to 0}} \mathcal{D}_\delta(\Lambda^*)=0$. 
The implication (ii)$\Rightarrow$(i) is immediate. Moreover, (i)$\Rightarrow$(ii) follows from
 $(iii) \Longrightarrow (ii)$, since the Khinchin action is free (Step 1).
The implication $(ii) \Longrightarrow (iii)$ is obvious. 
\end{proof}

\section{Obstruction by homomorphism into $\mathbb{R}$}
\label{sec:counting-functions}

Before proving Theorem~\ref{thm:subexponential-growth-obstruction}, we prove 
Theorem~\ref{thm:separated-prefix-concentration} which shows that the existence of a  homomorphism from $\mathbb{N}$ into $\mathbb{R}$ with some concentration along 
$\lambda_n$ is an obstruction for $(\lambda_n)$ to be $L^\infty$-Khinchin. 

\subsection{Proof of Theorem~\ref{thm:separated-prefix-concentration}: Obstruction by homomorphism to $\R$}

The proof of  Theorem~\ref{thm:separated-prefix-concentration} is based on Theorem~\ref{thm:homomorphism-criterion}.

\begin{proof}[Proof of Theorem~\ref{thm:separated-prefix-concentration}]
By Theorem~\ref{thm:homomorphism-criterion} (ii), it suffices to construct, {for each fixed integer $K>6$}, a measure-preserving multiplicative action $\gamma$ on some probability space
$(\mathbb{R}/T\mathbb{Z}, \mathfrak{m}_T)$ and a measurable set $E'$ such that
 $\|\one_{E'}\|_1=\mathfrak{m}_T(E')$ tends 
to $0$ as $K \to \infty$, whereas $\|\Lambda^{\gamma,*}\one_{E'}\|_1$ is bounded below by a positive
constant independent of $K$.  Here $T$ {will depend on} $K$ and $\mathfrak{m}_T$ is the normalized Haar measure on 
$\mathbb{R}/T\mathbb{Z}$.

The concentration assumption made on $\ell(\lambda_n)$ will first give a set $E$ with large maximal averages on a set $F$ satisfying
$\mathfrak{m}_T(E)/\mathfrak{m}_T(F)=3/K$.  Random translations will then enlarge $F$ to a set $F'$ of uniformly positive measure at a cost 
of only $O(K^{-1})$ for $E$.

Fix $K>6$.  By the hypothesis, we can choose $h>0$, $t_1,\ldots,t_K\in\mathbb R$, and $N_1,\ldots,N_K\in\mathbb N$ such that  $|t_i-t_j|\ge h \ \text{ for}\  i\ne j$ and
\begin{equation}\label{eq:concentration} 
         \frac1{N_i}\#\{n\le N_i:|\ell(\lambda_n)-t_i|\le h\}\ge\eta  \quad(1\le i\le K).
\end{equation}
Assume that $t_1<t_2<\cdots <t_K$ and take 
\[
       T>\max\bigl\{3h,\,t_K- t_1+h\bigr\}.
\]
The given homomorphism $\ell$ induces the measure-preserving $(\mathbb N,\cdot)$-action on the circle $\mathbb{T}_T:= \mathbb{R}/T\mathbb{Z}$ defined by 
\[
       \gamma^m(t)= t -\ell(m) \mod T,\qquad m\in\mathbb N.
\]
 Let $\rho_T:\mathbb R\to\mathbb T_T$ be the canonical projection, namely $\rho_T(x) = x \!\mod T$.  Introduce two sets in the circle $\mathbb{T}_T$ by projecting:
\[
     E=\rho_T([-h,2h)),\qquad       F=\rho_T\!\left(\bigcup_{i=1}^K[t_i,t_i+h)\right).
\]
Since the intervals $[t_i, t_i +h)$ on the line $\mathbb{R}$ which  define $F$ are pairwise disjoint, the choice of
$T$ prevents overlap after projection; hence, 
\begin{equation}\label{eq:E/F}
         \mathfrak{m}_T(E)=\frac{3h}{T},\qquad
         b:=\mathfrak{m}_T(F)=\frac{Kh}{T}\in(0,1),\qquad
         \frac{\mathfrak{m}_T(E)}{b}=\frac3K.
\end{equation}

Let $N_*:=\max_{1\le i\le K}N_i$.  For $y\in F$, there exist $i$ and $u\in[t_i,t_i+h)$ such that  $y=\rho_T(u)$.  For those $n\le N_i$ with
$|t_i - \ell(\lambda_n)|\le h$, we have $\gamma^{\lambda_n} y \in E$ because
\[
            -h\le t_i - \ell(\lambda_n)\le u-\ell(\lambda_n) < t_i+h -\ell(\lambda_n)\le 2h.
\]
This, together with the concentration hypothesis \eqref{eq:concentration}, yields

\[
           \forall \ y\in \rho_T([t_i,t_i+h)), \quad \Lambda_{N_i}^{\gamma}\one_E(y)
       = \frac1{N_i}\sum_{n=1}^{N_i}\one_E(\gamma^{\lambda_n}y) 
       \ge\eta.
\] 
It follows that the maximal function  $\Lambda_{N_*}^{\gamma,*}g  :=\max_{1\le N\le N_*}|\Lambda_N^\gamma g|$ satisfies
\begin{equation}\label{eq:finite-max-lower-bound}
            \Lambda_{N_*}^{\gamma,*}\one_E\ge\eta\,\one_F.
\end{equation}

We now use random translations to make the set of large maximal values have uniformly positive measure.  Put $M=\lceil 1/b\rceil$ and choose $M$
independent $\mathfrak{m}_T$-uniform random points $\xi_1,\ldots,\xi_M$ in $\mathbb T_T$.  For every $t \in\mathbb T_T$,
\[
     \mathbb P\left(\omega: t \notin\bigcup_{j=1}^M(F+\xi_j(\omega))\right)
 =(1-b)^M\le e^{-Mb}\le e^{-1}.
\]
Fubini's theorem shows that  there is a deterministic choice of  
translates $s_j =\xi_j(\omega^*)$ for some $\omega^*$ for which
\[
      F':=\bigcup_{j=1}^M(F+s_j),  \qquad \mathfrak{m}_T(F')\ge1-e^{-1}.
\]
Let $E':=\bigcup_{j=1}^M(E+s_j)$.  By the fact $M\le b^{-1}+1<2b^{-1}$ and the third equality in \eqref{eq:E/F}, {we obtain}
\begin{equation}\label{eq:E'}
          \|\one_{E'}\|_1= \mathfrak{m}_T(E')\le M \mathfrak{m}_T(E)<\frac{6}{K}.
\end{equation}
On the other hand, for each $y\in F'$, there exists  $j$ such that  $y-s_j\in F$.   
As $\one_{E'}\ge \one_{E+s_j}$, we get
\[
      \Lambda_{N_*}^{\gamma,*}\one_{E'}(y)
   \ge\Lambda_{N_*}^{\gamma,*}\one_{E+s_j}(y)
      =\Lambda_{N_*}^{\gamma,*}\one_E(y-s_j)
    \ge\eta
\]
where the equality holds because the translations  commute  with $\gamma$, and where the last inequality is ensured by
~\eqref{eq:finite-max-lower-bound}.  Integrating the last inequality over $F'$ leads to
\[
      \|\Lambda^{\gamma,*}\one_{E'}\|_1
   \ge \eta(1-e^{-1}).
\]
This, together with \eqref{eq:E'}, implies 
\[
     \sup_{\gamma}\mathcal D_{6/K}(\Lambda^{\gamma,*})
 \ge\eta(1-e^{-1})\qquad(K>6).
\] 
Letting $K\to\infty$ contradicts the uniform continuity condition in Theorem~1.1(ii).  Hence $\Lambda$ is not $L^\infty$-Khinchin.
\end{proof}

\subsection{Proof of Theorem~\ref{thm:subexponential-growth-obstruction}: Necessary growth of $L^\infty$-Khinchin sequences}

\begin{proof}[Proof of Theorem \ref{thm:subexponential-growth-obstruction}]
We apply Theorem~\ref{thm:separated-prefix-concentration} with
$\ell(n)=\log n$ and $\eta=1/2$. For $i\ge0$, put
$b_i=\log\lambda_{2^i}$. For $i\ge1$, let $t_i$ be the midpoint of the interval
$J_i=(b_{i-1},b_i]$. Since  $\lambda_n$ is strictly increasing, for $2^{i-1}<n\le2^i$ we have $\log\lambda_n\in J_i$ so that 
\begin{equation}\label{eq:count}
 \frac{1}{2^i}\#\{n\le2^i:\log\lambda_n\in J_i\}\ge\frac12.
\end{equation}
Fix an arbitrary integer $K\ge1$. In order to apply Theorem~\ref{thm:separated-prefix-concentration},  it suffices to find $h>0$ and $K$ intervals $J_i$ of
length at most $2h$ whose midpoints are pairwise at least $h$ apart.

We distinguish two cases.

If infinitely many $J_i$ have length less than $1$, we take $h=1$.
Since $t_i\to\infty$, we can successively choose $K$ intervals from these intervals $J_i$ such that the midpoints of the chosen intervals have distances bigger than $1$.

Suppose now that only finitely many intervals $J_i$ have length  $|J_i|$ less than $1$.
Let $C=\#\{i\ge1:|J_i|<1\}$ be the number of these intervals. By the hypothesis, 
we can choose $n_r\to\infty$ such  that
\begin{equation}\label{eq:hypothesis2}
        \log\log\lambda_{n_r}=o(\log n_r).
\end{equation}
Let
$m_r=\lfloor\log_2 n_r\rfloor$ be the integer such that $2^{m_r}\le n_r<2^{m_r +1}$.
We have $\log n_r\sim m_r\log 2$. Then from \eqref{eq:hypothesis2} we get
\begin{equation}\label{eq:b-m}
 \log b_{m_r}
 =\log\log\lambda_{2^{m_r}}
 \le\log\log\lambda_{n_r}=o(m_r).
\end{equation}

There are at least $m_r-C$ intervals among $J_1,\ldots,J_{m_r}$ whose lengths fall
into the interval $[1,b_{m_r}]$. Cover the interval $[1, b_{m_r}]$ by disjoint dyadic intervals
 $[2^k,2^{k+1})$ 
($0\le k\le  \lfloor\log_2 b_{m_r}\rfloor$). 
The pigeonhole principle gives
\begin{equation}\label{eq:pigeon}
    \max_{0\le k\le \lfloor\log_2 b_{m_r}\rfloor}\#\{1\le i\le m_r:2^k\le |J_i|<2^{k+1}\}
          \ \ge\ \frac{m_r-C}{1+\lfloor\log_2 b_{m_r}\rfloor}.
\end{equation}
Since
$1+\lfloor\log_2 b_{m_r}\rfloor=o(m_r)$ according to \eqref{eq:b-m}, the 
right-hand side of \eqref{eq:pigeon} tends to infinity as $r\to \infty$.  
Therefore, for a large $r$, we can choose  $0\le k_0\le \lfloor\log_2 b_{m_r}\rfloor$ such that for $h=2^{k_0}$ we have
$$
 \#\{1\le i\le m_r: h\le |J_i|<2 h\}
          \ \ge\ K.
$$
That means, we can select $K$ (disjoint)  intervals $J_i$ with lengths between $h$ and $2h$
with $h\le |J_i|<2h$. 

{In either case, each selected interval satisfies
$J_j\subseteq[t_j-h,t_j+h]$. Taking $N=2^j$ for each selected interval, \eqref{eq:count} verifies
Theorem~\ref{thm:separated-prefix-concentration} with $\eta=1/2$ }
\end{proof}

\color{black}

The following are some direct applications of  Theorem \ref{thm:subexponential-growth-obstruction}.
\begin{example}
 The following sequences are not
  $L^\infty$-Khinchin:
  \[
    (n)_{n\ge1},\qquad
    \bigl(n^{\lfloor\log n\rfloor}\bigr)_{n\ge3},\qquad
    \bigl(n^{\lfloor\log\log n\rfloor}\bigr)_{n\ge16}.
  \]
\end{example}

\begin{example}
Let $p_n$ be the $n$th prime. The prime number theorem
  gives $p_n\sim n\log n$, so $\log\log p_n/\log n\to0$.
  Theorem~\ref*{thm:subexponential-growth-obstruction}
thus recovers Marstrand's conclusion
  \cite[Theorem~5.3 and Corollary~3.3]{Marstrand} that the increasing sequence of all primes is not
  $L^\infty$-Khinchin.
\end{example}

\subsection[A primitive substitution producing a non-\texorpdfstring{$L^\infty$}{L-infinity}-Khinchin sequence]
{A primitive substitution producing a non-$L^\infty$-Khinchin sequence}
In~\cite[Proposition~6.1]{FanFanLiLp2026},  we proved that 
any $C$-balanced primitive {sequence} produces  an $L^1$-Khinchin sequence. 
In the opposite direction, 
{using} Theorem~\ref{thm:separated-prefix-concentration}, we are now going  to construct
a primitive substitution whose fixed point produces a non-$L^\infty$-Khinchin
sequence.

\color{black}
\begin{corollary}\label{counter-example-primitive}
Let $w=(w_n)_{n\ge1}$ be the fixed point of the primitive substitution
\[
 \sigma(0)=010,\qquad \sigma(1)=011.
\]
Then, for coprime integers $q_0,q_1\ge2$, the sequence
\[
 \Lambda(w)=(q_{w_1}\cdots q_{w_n})_{n\ge1}
\]
is not an $L^\infty$-Khinchin sequence.
\end{corollary}

\begin{proof}
  Write $\lambda_n=q_{w_1}\cdots q_{w_n}$. To apply Theorem~\ref{thm:separated-prefix-concentration}, we choose a
homomorphism $\ell:\mathbb{N}\to \mathbb{R}$ such that $\ell(\lambda_n)$ equals 
{the} difference of the number of zeros
and the number of ones in $w_1\cdots w_n$. 
Denote this by 
\[
  D_0=0; \qquad D_n=\sum_{j=1}^n (1-w_j)-\sum_{j=1}^n w_j =\sum_{j=1}^n(1-2w_j) \quad\text{for} \ n \ge 1.
\]

Choose primes $p_0\mid q_0$ and $p_1\mid q_1$. By coprimality, 
the homomorphism
\[
 \ell(m):=\frac{v_{p_0}(m)}{v_{p_0}(q_0)}
        -\frac{v_{p_1}(m)}{v_{p_1}(q_1)} \quad (m \in \mathbb{N})
\]
is well defined and 
satisfies $\ell(\lambda_n)=D_n$.

Since $\sigma(a)=01a$, we have
\[
 D_{3m}=D_m,\qquad
 D_{3m+1}=D_m+1,\qquad
 D_{3m+2}=D_m.
\]
Together with $D_0=0$, these identities show that $D_n$ is the
number of digits equal to $1$ in the ternary expansion of $n$.

For $n$ chosen uniformly from $\{0,\ldots,3^r-1\}$, its $r$
ternary digits are independent and uniform on $\{0,1,2\}$.
Thus $D_n$ has mean $r/3$ and variance $2r/9$, giving
\[
 \frac1{3^r-1}\sum_{n=1}^{3^r-1}
       \left(D_n-\frac r3\right)^2
 \le \frac{3^r}{3^r-1}\frac{2r}{9}
 \le\frac r3
 \qquad(r\ge1).
\]

To obtain $K$ separated centers with a common width, fix $K\ge1$
and set
\[
 L=6K,\qquad
 N_i=3^{iL}-1,\qquad
 t_i=\frac{iL}{3},\qquad h=\frac L3
 \qquad(1\le i\le K).
\]
Then $|t_i-t_j|\ge h$ for $i\ne j$. Applying the preceding estimate
with $r=iL$, Markov's inequality yields
\[
 \frac1{N_i}\#\{n\le N_i:|\ell(\lambda_n)-t_i|>h\}
 \le\frac{iL/3}{h^2}
 =\frac{3i}{L}\le\frac12.
\]
Hence Theorem~\ref{thm:separated-prefix-concentration} applies
with $\eta=1/2$.
\end{proof}

\color{black}

\section{Obstruction by homomorphism into $\mathbb{T}$}
\label{sec:examples-homomorphism-obstruction}

Using his entropy method, Bourgain \cite[Proposition 10]{Bourgain1988} proved the following important result:  for any 
sequence of real numbers $a_k\to a$ with $a_k \not =a$ for  each $k$ (it was stated for $a=0$), there exists a bounded function
$f\in L^\infty(\T)$ such that
\begin{align*}
\frac1n\sum_{k=1}^n f(x+a_k) \ \ \text{diverges on a set of positive measure}.
\end{align*}
 Kostyukovsky--Olevskii
\cite{KostyukovskyOlevskii2005} extended the result by using  discrete
probability measures.
Karagulyan \cite{Karagulyan2011} further strengthened it to strong sweeping out.
We shall use Theorem~\ref{thm:homomorphism-criterion}, together with Karagulyan's result, 
to construct obstructions to the $L^\infty$-Khinchin property.

Karagulyan is concerned with convolution operators. Let us recall
the definition of the convolution of a measure $\mu$ with a function $f$:
\begin{align*}
\mu*f(x):=\int_\T f(x-t)\,d\mu(t).
\end{align*}

\begin{lemma}[Karagulyan, \cite{Karagulyan2011}, Theorem 2]\label{Karagulyan}
Let $(\nu_n)_{n\ge1}$ be discrete probability measures on $\T$.
Suppose that, for some $a\in\T$,
\begin{align*}
\nu_n(\{a\})=0\quad( \forall \ n\ge1);
\qquad
\lim_{n\to\infty}\nu_n(B(a,r))=1\quad(\forall \ 0<r\le1/2).
\end{align*}
Then, 
for every $\varepsilon>0$, there is a measurable set $E\subset\T$
with $\m(E)<\varepsilon$ such that
\begin{align*}
\limsup_{n\to\infty}\nu_n*\one_E(x)=1,
\qquad
\liminf_{n\to\infty}\nu_n*\one_E(x)=0
\end{align*}
for $\m$-almost every $x$.
\end{lemma}

\subsection{Proof of Theorem \ref{Bourgain-obstruction}: Obstruction by homomorphism into $\mathbb{T}$}

\begin{proof}
Consider the discrete probability measures describing the averages of $(\pi(\lambda_n))$:
\[
\mu_n:=\frac1n\sum_{k=1}^n\delta_{\pi(\lambda_k)},
\qquad n\ge1.
\]
Denote \(
B^\circ(a,r):=B(a,r)\setminus\{a\}
\) for $0<r<1/2$, which is a punctured interval.
The assumption \eqref{eq:K-cond} in Theorem~\ref{Bourgain-obstruction} implies that
\begin{equation}\label{notation:kappa}
\kappa:=\lim_{{r\to 0}}\limsup_{n\to\infty}\mu_n(B^\circ(a,r)) \ge \overline{d}(I)>0.
\end{equation}

Choose $r_j\to0$ with $0<r_j<1/2$, so that 
$$b_j:=\limsup_{n\to\infty}\mu_n(B^\circ(a,r_j)) \to\kappa.
$$ By the definition of $b_j$, we can recursively choose
integers $n_j>n_{j-1}$  so that $c_j:=\mu_{n_j}(B^\circ(a,r_j))$ satisfies
\begin{align*}
|c_j-b_j|<1/j.
\end{align*}
Observe that $c_j\to\kappa$. Since $\kappa>0$,  $c_j>0$ for large $j$, say $j\ge j_0$. Define the probability measures
\[
\nu_j
:=
\frac{\mu_{n_j}|_{B^\circ(a,r_j)}}{c_j} \quad \text{for } j\ge j_0.
\]
The sequence of measures $(\nu_j)_{j\ge j_0}$ clearly satisfies the condition of Lemma~\ref{Karagulyan}, which implies that, 
for  $\varepsilon>0$, there exists a measurable set
$E\subset\T$ with $\m(E)<\varepsilon$ such that
\begin{equation}\label{eq:K2}
\limsup_{j\to\infty}\nu_j*\one_E(x)=1
\qquad\text{for $\m$-almost every }x.
\end{equation}
Since $\one_E\ge0$ and
$\mu_{n_j}|_{B^\circ(a,r_j)}=c_j\nu_j$, we have
\[
\mu_{n_j}*\one_E(x)\ge c_j\,\nu_j*\one_E(x).
\]
Then \eqref{eq:K2} implies
\begin{align}
\limsup_{j\to\infty}\mu_{n_j}*\one_E(x)
&\ge
\limsup_{j\to\infty}c_j\,\nu_j*\one_E(x)\ge \lim_{j\to \infty}c_j
=\kappa
\qquad\text{for $\m$-a. e.  }x.
\label{eq:K3}
\end{align}

This  can be interpreted as a sweeping out result for the measure-preserving {$(\N,\cdot)$-action} on $\T$ defined by $\pi(n)$-translation
\[
{\gamma^n} x:=x-\pi(n)\pmod1.
\]
Indeed, 
since 
$\Lambda_n^{{\gamma}} f=\mu_n*f$, \eqref{eq:K3} implies 
\[
\Lambda^{{\gamma},*}\one_E(x)\ge\kappa
\qquad\text{for $\m$-a. e. }x.
\]
This implies that the maximal operator $\Lambda^{{\gamma},*}$ is not continuous at $0$ on $B_\infty(\T)$, because 
$\m(E)=\|\one_E\|_1$ can be {arbitrarily small}.
Therefore the condition (ii) in Theorem~\ref{thm:homomorphism-criterion} is violated,   so
$\Lambda$ is not an $L^\infty$-Khinchin sequence by Theorem~\ref{thm:homomorphism-criterion}.
\end{proof}

\subsection[When the homomorphic images have finitely many subsequential limits]
{When the homomorphic images
\texorpdfstring{$(\pi(\lambda_n))_{n\ge1}$}{}
have finitely many subsequential limits}

We are going to show that a sequence of integers is not  $L^\infty$-Khinchin if  
its image under a homomorphism has finitely many limit points and eventually avoids its limit points.
We recall  that 
 the limit of a convergent subsequence of a sequence is  called a {\em limit point} of the sequence.

\begin{corollary}\label{cor-of-Bourgain-obstruction}
Let $\Lambda=(\lambda_n)_{n\ge1}$ be a sequence of distinct positive
integers. 
Suppose that there exists a homomorphism $\pi:(\N,\cdot)\to(\T,+)$ such that the sequence $(\pi(\lambda_n))$ has exactly $m$ subsequential limits $a_1,\ldots,a_m$ and that
\begin{align*}
\pi(\lambda_n)\notin\{a_1,\ldots,a_m\}
\qquad\text{for all sufficiently large }n.
\end{align*}
Then neither $\Lambda$ nor any rearrangement of $\Lambda$ is an $L^\infty$-Khinchin sequence.
\end{corollary}

This is a corollary of Theorem~\ref{Bourgain-obstruction}. To prove it we  use the following elementary fact.

\begin{lemma}\label{lm:upper-density}
  Let $(x_n)_{n \ge 1}$ be a sequence in $\T$.
  Assume that $(x_n)$ has exactly $m$ limit points $a_1, \ldots, a_m \in \T$.
  Then there exist a subset $I$ of $\N$ and $a_j$ such that
  \begin{align*}
    \lim_{\substack{n \to \infty \\ n \in I}} x_n= a_j, \qquad \overline{d}(I)  \ge \frac{1}{m}.
  \end{align*}
\end{lemma}

\begin{proof}
Choose open neighborhoods $U_1,\ldots,U_m$ of $a_1,\ldots,a_m$, respectively, with pairwise disjoint closures.
Since $a_1, \cdots, a_m$ are the only limit points,  there are only finitely many terms $x_n$ outside the union $\bigcup_{j=1}^m U_j$. 
So,  the sets
$I_j:=\{n:x_n\in U_j\}$ form a partition of all but finitely many positive integers.
Consequently,  at least one $I_j$ has upper density at least $1/m$.
For such a $j$, the subsequence $ (x_{i})_{i\in I_j}$ converges to $a_j$. 
Otherwise, infinitely many of its terms would lie in the compact set
$\overline{U}_j\setminus B(a_j,\varepsilon)$ for some $\varepsilon>0$,
which would produce another
 limit point in $\overline{U}_j$, a contradiction. 
\end{proof}

\begin{proof}[Proof of Corollary~\ref{cor-of-Bourgain-obstruction}]
Applying Lemma~\ref{lm:upper-density} to $x_n=\pi(\lambda_n)$
yields  a subset  $I\subset\mathbb{N}$ with $\overline d(I)\ge1/m$ such that 
$\pi(\lambda_n)\to a_j$ for some $j$ as $n \in I$ tends to $\infty$. 
Together with the eventual-avoidance assumption, this implies that
for each $r>0$, we have
$\pi(\lambda_n)\in B(a_j,r)\setminus\{a_j\}$ for all but
finitely many $n\in I$.

After deleting finitely many indices from $I$, condition \mbox{\eqref{eq:K-cond}} holds, 
so Theorem~\mbox{\ref{Bourgain-obstruction}} applies.

A rearrangement of $\Lambda$ does not change the set of subsequential limits, nor the property of eventually avoiding this  set of limits. 
Hence every rearrangement of $\Lambda$ satisfies the hypotheses of Corollary ~\ref{cor-of-Bourgain-obstruction}.
Applying the preceding argument to the rearranged sequence, we conclude that every rearrangement of $\Lambda$ is also non-$L^\infty$-Khinchin.
\end{proof}

We have the following class of non-$L^\infty$-Khinchin sequences. 

\begin{corollary}\label{prop:leading-part}
Let $\Lambda=(\lambda_n)_{n\ge1}$ be a sequence of distinct positive integers of
{
the form $\lambda_n=q^{a(n)}r_n$, where  $a(n) \in \Z$, $q>1$ and $r_n>0$.}
Suppose that $(\log_qr_n\bmod1)$ has a finite set $A$ of
limit points and that
\begin{align*}
\log_qr_n\bmod1\notin A
\qquad\text{for all sufficiently large }n.
\end{align*}
Then $\Lambda$, as well as any rearrangement,  is not $L^\infty$-Khinchin.
\end{corollary}

\begin{proof}
We consider the homomorphism $\pi(n):=\log_qn\bmod1$ from $\mathbb{N}$ to $\mathbb{T}$. Then
$$\pi(\lambda_n)=\log_qr_n\bmod1.
$$
Hence, Corollary~\ref{cor-of-Bourgain-obstruction} ensures the conclusion.
\end{proof}

More concretely, we have the following classes of examples.

\begin{example} A sequence of distinct positive integers of the form $q^{a(n)} +b(n)$ with
a positive number $q>1$, a sequence of integers $a(n)$  and a  sequence of positive numbers $b(n)$ 
 such that $b(n) =o(q^{a(n)})$
is not $L^\infty$-Khinchin. In particular,
$(3^n+2^n)_{n\ge 1}$ is not $L^\infty$-Khinchin.
\end{example}

Indeed,  $q^{a(n)} +b(n)$ is of the form $q^{a(n)}r_n$ with  $\log_q r_n = \log_q (1+b(n)q^{-a(n)}) \to 0$ but $\log_q r_n \not =0$. 

\begin{example}\label{exam:leading-part} Let $q\ge2$ be an integer. Then
$(q^n+1)_{n\ge 1}$ is not $L^\infty$-Khinchin. 
\end{example}


\begin{example}\label{prop:union-of-non-K}
For any finite nonempty set $I\subset\N_0$, every enumeration of
$
\{2^t3^m+1:t\in I,\ m\ge0\}
$
is not $L^\infty$-Khinchin.
\end{example}

Indeed, if we
consider the homomorphism $\pi(n)=\log_3n\bmod1$, we have
\begin{align*}
\pi (2^t3^m+1)
=t\log_3 2+\log_3(1+2^{-t}3^{-m})\pmod1.
\end{align*}
Hence the set of limit points of the image sequence is
 the finite set 
$\{t\log_3 2\bmod1: t \in I\}$.
Since $\log_3 2$ is an irrational number, $t\log_3 2\bmod1 (t \in I)$ are pairwise
distinct. 
The avoidance condition is satisfied and
Corollary~\ref{cor-of-Bourgain-obstruction} applies.

Let  $(s_n)_{n\ge1}$ be the increasing enumeration
of the Furstenberg set $S(2,3):=\{2^n3^m:n,m\ge0\}$. The translated sequence $(s_n+1)_{n\ge1}$ is 
not $L^\infty$-Khinchin. More generally, any sequence of distinct numbers from $S(2,3)+1$ is 
not $L^\infty$-Khinchin, see {Example \ref{shift-Furstenberg}} below.

\begin{corollary}\label{cor:prime-subsequence-not-khinchin}
Let $(p_n)_{n\ge1}$ be any sequence of distinct primes. Neither of the
sequences
$
(p_n)_{n\ge1}$ and $\left(\prod_{j=1}^n p_j\right)_{n\ge1}
$
is $L^\infty$-Khinchin.
\end{corollary}

\begin{proof} Let $\mathcal{P}$ be the set of all prime numbers. We check that Corollary~\ref{cor-of-Bourgain-obstruction} applies  in both cases.
Recall that any homomorphism  $\pi: \mathbb{N}\to \mathbb{T}$ is uniquely determined by its values on primes, and these values can be chosen arbitrarily.

 To study  the sequence $(p_n)_{n\ge 1}$, we consider the homomorphism $\pi_1$ defined by 
$$
\pi_1(p_n)=\frac{1}{n+1} , \quad \pi_1(p) =0 \ \text{for} \ p \in \mathcal{P}\setminus \{p_n: n \ge 1\}.
$$  
Clearly, $\pi_1(p_n)\to 0$  as $n\to \infty$ and $\pi_1(p_n) \not=0$ for all $n\ge 1$. 

To study the sequence $(p_1p_2\cdots p_n)_{n\ge 1}$, we consider the homomorphism $\pi_2$ defined by
\begin{align*}
\pi_2(p_1)=\frac{1}{2},\ \
\pi_2(p_n)=\frac1{n+1}-\frac1n\ \ (n\ge2),  \ \ \pi_2(p) =0 \ \text{for} \ p \in \mathcal{P}\setminus\{p_n: n \ge 1\}.
\end{align*}
Then
$\pi_2(p_1\cdots p_n)=1/(n+1) \to 0$.
\end{proof}

In Corollary~\ref{cor:prime-subsequence-not-khinchin}, we only need to assume that $(p_n)_{n\ge1}$ is a sequence of pairwise coprime positive integers with ${p_n \ge 2}$.
The proof is the same, but we need Lemma~\ref{lem:homom} for ensuring the existence of a suitable homomorphism.

\section{Obstruction by homomorphism into $\mathbb{Q}_p$}

We derive the $p$-adic obstruction from the circle homomorphism
criterion (Theorem~\ref{Bourgain-obstruction}). 

\subsection{Proof of Theorem~\ref{thm:padic-homomorphism-obstruction}: Obstruction by homomorphism into $\mathbb{Q}_p$}
We shall deduce Theorem \ref{thm:padic-homomorphism-obstruction}  from
Theorem \ref{Bourgain-obstruction}.
The key point is  the following lemma which constructs specific group homomorphisms from $\mathbb{Q}_p$ into $\mathbb{T}$.
Such group homomorphisms belong to the Bohr compactification of {$\widehat{\mathbb{Q}_p}$}. The proof of the lemma will be  based on Baire's theorem.

\begin{lemma}\label{lem:padic-character-selection}
Let $(a_n)_{n \ge 1}$ be a sequence in $\Q_p$.
Suppose there exists a set $I\subseteq\N$ of positive upper density such that
$a_n \neq 0$  for all $ \ n \in I$ and 
\begin{equation*}
    \lim_{\substack{n \to \infty\\ n \in I}} a_n = 0.
\end{equation*}
Then there exists a group homomorphism $\chi:(\Q_p,+)\to(\T,+)$ and a set
$I'\subseteq\N$ of positive upper density such that $\chi(a_n) \neq 0$ for all $ n \in I'$ and 
\[
 \lim_{\substack{n\to\infty\\ n\in I'}}\chi(a_n)=0. 
\]
\end{lemma}


\begin{proof} The group homomorphism $\chi$ that we are looking for is not assumed continuous. So, if it does exist, it belongs to the Bohr compactification group $b\widehat{\mathbb{Q}_p}$  of $\widehat{\mathbb{Q}_p}$.
Recall that 
\[
 b\widehat{\Q_p}=\{ \chi \in \T^{\mathbb{Q}_p}: \chi(h_1+ h_2) =\chi(h_1) + \chi(h_2), \ \forall \ h_1, h_2 \in \Q_p \}.
\] 
Here $\T^{\Q_p}$ is the space of all functions defined on $\mathbb{Q}_p$ taking values in $\mathbb{T}$. It is nothing but a product space of $\mathbb{T}$'s, a compact Hausdorff group with the topology of
 pointwise convergence. Clearly,  $b\widehat{\Q_p}$ is a closed subgroup of
$\T^{\Q_p}$ and is therefore compact Hausdorff.
We will use  the Baire category
theorem to confirm the existence of the desired group homomorphism $\chi\in b\widehat{\Q_p}$.

We shall use  two   known facts. The first one is that $ \widehat{\Q_p} $ is dense 
in $b\widehat{\Q_p}$ (\cite[Theorem~1.8.2]{Rudin1990}).
The second one  is that  for each $h\in\Q_p\setminus\{0\}$, the closed subgroup
\[
 Z_h=\{\chi\in b\widehat{\Q_p}:\chi(h)=0\}
\]
has empty interior,  by subgroup--quotient duality
\cite[Theorem~2.1.2]{Rudin1990}, which identifies
$b\widehat{\Q_p}/Z_h$ with the circle group $\T$.

The assumption made on $(a_n)$ implies that
\begin{equation}\label{cond:padic-concentration*}
 \kappa:=\lim_{k\to\infty}\limsup_{N\to\infty}
 \frac1N\#\bigl\{1\le n\le N:
       a_n \in
       p^k\mathbb Z_p\setminus\{0\}\bigr\}>0.
\end{equation}
Here $p^k\mathbb Z_p\setminus\{0\}$ is a punctured ball with center $0 \in \Q_p$.

Fix $0<\eta<\kappa$. For $0<\varepsilon<1/2$ and $K\in\N$, define
\[
 \mathcal U_{\varepsilon,K}
 =\bigcup_{N\ge K}
 \left\{\chi\in b\widehat{\Q_p}:
 \frac1N\#\{n\le N:0<\|\chi(a_n)\|_\T<\varepsilon\}>\eta
 \right\}.
\]
We claim that each $ \mathcal U_{\varepsilon,K}$ is open and dense.
The openness follows from the expression
\[
 \mathcal U_{\varepsilon,K}
 =
 \bigcup_{N\ge K}
 \bigcup_{\substack{J\subseteq\{1,\ldots,N\}\\ |J|>\eta N}}
 \ \bigcap_{n\in J}
 \left\{\chi\in b\widehat{\Q_p}:
       0<\|\chi(a_n)\|_\T<\varepsilon\right\}.
\]
Let us  prove the density.  Let $W\subset b\widehat{\Q_p}$ be a nonempty open
set. By the first fact mentioned above, $W$ contains some character $\chi_\xi \in \widehat{\mathbb{Q}_p}$ 
defined by $\chi_{\xi}(x)=\psi_p(\xi x)$ 
where $\xi \in \Q_p$ and $\psi_p$ is the standard character 
on $(\Q_p,+)$ with kernel $\Z_p$.

Take a large integer $k\ge1$
such that $\xi p^k\mathbb Z_p\subset\mathbb Z_p$. Since
$\eta<\kappa$, by \eqref{cond:padic-concentration*}, there exists an integer $N\ge K$ for which
\[
 J:=\{1\le n\le N:a_n\in p^k\mathbb Z_p\setminus\{0\}\}
 \qquad\text{satisfies}\qquad |J|>\eta N.
\]
For each $n\in J$, we have $\xi a_n \in  \xi p^k \mathbb{Z}_p \subset \mathbb{Z}_p$ so that $\chi_\xi(a_n) =0$. Consequently,
$\chi_\xi \in W'$ where 
\[
 W':=W\cap\bigcap_{n\in J}
       \{\chi\in b\widehat{\Q_p}:\|\chi(a_n)\|_\T<\varepsilon\}
\]
which is open. 
 For every $n\in J$,  $a_n\ne0$, hence by the second fact mentioned above,
$b\widehat{\Q_p}\setminus Z_{a_n}$ is open and dense. 
Since $J$ is
finite, the intersection of these sets $b\widehat{\Q_p}\setminus Z_{a_n}$ is dense and therefore
meets the nonempty open set $W'$. Thus
\[
 W'\cap\bigcap_{n\in J}
       \bigl(b\widehat{\Q_p}\setminus Z_{a_n}\bigr)
 \ne\varnothing, \ \ \ \text{namely} \ \ \  W'\setminus\bigcup_{n\in J}Z_{a_n}  \ne\varnothing.
\] 
Every character in this set satisfies
$0<\|\chi(a_n)\|_\T<\varepsilon$ for all $n\in J$. Since
$|J|>\eta N$, each such character belongs to
$W\cap\mathcal U_{\varepsilon,K}$. This proves density.

Take a sequence $(\varepsilon_j)_{j \ge 1}$ with $0<\varepsilon_j<1/2$ such that $\varepsilon_j\to0$.
Since $b\widehat{\Q_p}$ is compact Hausdorff, we can apply the Baire category
theorem to the open dense sets $\mathcal U_{\varepsilon_j,K}$ $(j,K \in \N)$ to give a character
\[
 \chi\in\bigcap_{j\ge1}\bigcap_{K\ge1}\mathcal U_{\varepsilon_j,K}.
\]
For every $j$ there are arbitrarily large $N$ such that
\[
 \#\{1\le n\le N:0<\|\chi(a_n)\|_\T<\varepsilon_j\}>\eta N.
\]
Set $N_0=0$. Recursively choose such $N_j$ with $N_j>N_{j-1}$
and $N_{j-1}<\eta N_j/2$. Then put
\[
 I_j=\{N_{j-1}<n\le N_j:
                  0<\|\chi(a_n)\|_\T<\varepsilon_j\},
 \qquad I'=\bigcup_{j\ge1}I_j.
\]
Then
\[
 \#(I'\cap[1,N_j])
 \ge |I_j|>\eta N_j-N_{j-1}>\frac{\eta N_j}{2}.
\]
Hence $I'$ has upper density at least $\eta/2>0$. Moreover,
$\chi(a_n)\ne0$ on $I'$ and 
$\chi(a_n)\to0$ along $I'$.
\end{proof}

\begin{proof}[Proof of Theorem~\ref{thm:padic-homomorphism-obstruction}]
Put
\(
 a_n=\pi_p(\lambda_n)-a \in \Q_p, n\ge1.
\)
Then the hypothesis made in Theorem~\ref{thm:padic-homomorphism-obstruction} 
can be stated as 
\begin{equation*}
    \lim_{\substack{n \to \infty\\ n \in I}} a_n = 0, \qquad a_n \neq 0 \ \  \text{for all} \ \ n \in I.
\end{equation*}
By Lemma~\ref{lem:padic-character-selection}, there exist a group
homomorphism $\chi:(\Q_p,+) \to(\T,+)$ and a set $I'\subset\N$ of positive upper
density such that
\[
 \chi(a_n)\to0\quad(n\to\infty,\ n\in I'),
 \qquad \chi(a_n)\ne0\quad(n\in I').
\]
Then the composed semigroup homomorphism
\(
 \pi:=\chi\circ\pi_p:(\N,\cdot)\longrightarrow(\T,+)
 \) verifies
\[
 \pi(\lambda_n)-\chi(a)
 =\chi\bigl(\pi_p(\lambda_n)-a\bigr)  \to 0 \quad(n\to\infty,\ n\in I'), \qquad  \pi(\lambda_n)\ne \chi(a) \quad(n\in I').
\]
Then we can apply Theorem~\ref{Bourgain-obstruction} to conclude  that $\Lambda$ is not an
$L^\infty$-Khinchin sequence.
\end{proof}

Theorem~\ref{thm:padic-homomorphism-obstruction} can be generalized as follows.  
This generalized form is not used in this paper, and its proof is omitted.
\begin{theorem}
\label{thm:G-homomorphism-obstruction}
Let $\Lambda=(\lambda_n)_{n\ge1}$ be a sequence
of distinct positive integers. Let $(G,+)$ be a  locally compact Hausdorff
 abelian group such that some neighborhood of the neutral element contains no nontrivial finite subgroup.
Suppose that there exist a semigroup homomorphism $\pi_G: (\N, \cdot) \longrightarrow (G, +)$, a point 
$a \in G$ and a subset $I (\subset \N)$ of positive upper density such that 
\begin{equation}\label{eq:G-concentration}
    \lim_{\substack{n \to \infty\\ n \in I}} \pi_G(\lambda_n) = a, \qquad \pi_G(\lambda_n) \neq a \ \text{for every} \ n \in I.
\end{equation}
Then $\Lambda$ is not an $L^\infty$-Khinchin sequence.
\end{theorem}

However, we can prove the 
special case $G=\mathbb{Q}_p^\times$ of  Theorem~\ref{thm:G-homomorphism-obstruction}
as a consequence of Theorem~\ref{thm:padic-homomorphism-obstruction}, using the fact that the logarithm is a homomorphism from
$\mathbb{Q}_p^\times$ into $\mathbb{Q}_p$.

\begin{theorem}\label{multiplicative-p-obs}
  Let $\Lambda=(\lambda_n)_{n \ge 1}$ be a sequence of 
  distinct positive integers. 
  Suppose there exist a homomorphism $\pi: (\N, \cdot) \to (\Q_p^\times, \cdot)$ 
  , $c \in \Q_p^\times$ and a set $I \subset \N$ of positive upper density
  such that 
 \begin{align*}
   \lim_{\substack{n \to \infty \\ n \in I} }\pi(\lambda_n) = c  \quad \text{in} \ \Q_p, \qquad \pi(\lambda_n) \neq c \ \  \text{for} \ \ n \in I.
 \end{align*}
 Then $\Lambda$ is not $L^\infty$-Khinchin.
\end{theorem}

\begin{proof}
  Recall that the standard \(p\)-adic logarithm \(
\log_p:\mathbb Q_p^\times\longrightarrow\mathbb Q_p,
\) normalized by \(\log_p(p)=0\),  is a continuous group homomorphism.
Moreover, the restriction of \(\log_p\) to
\[
U_p=
\begin{cases}
1+p\mathbb Z_p,&p\ge3,\\
1+4\mathbb Z_2,&p=2,
\end{cases}
\]
is injective. 
Consider the semigroup homomorphism from \((\mathbb N,\cdot)\) to \((\mathbb Q_p,+)\) defined by 
\[
\pi_p(n)=\log_p(\pi(n)),\qquad n\ge1,
\]
which is the composition of $\pi$ and $\log_p$.  By the hypothesis on $\pi$ and the continuity of \(\log_p\),
\[
\lim_{\substack{n\to\infty\\n\in I}}\pi_p(\lambda_n)=a, \quad \text{where} \ \ 
a:=\log_p(c). 
\]

To apply Theorem~\ref{thm:padic-homomorphism-obstruction} to conclude, 
it remains to show that $\pi_p(\lambda_n)\not=a$ for sufficiently large $n \in I$. 
As 
\(
\pi(\lambda_n)\) tends to $ c$ when $n\in I$ tends to infinity and $c\not=0$, 
\(
\frac{\pi(\lambda_n)}{c}\in U_p
\) for large $n\in I$. On the other hand, it is assumed that 
\(\pi(\lambda_n)\ne c\), hence $\frac{\pi(\lambda_n)}c\ne1$,  for large $n\in I$. 
By injectivity of \(\log_p\) on \(U_p\),
\[
\begin{aligned}
\pi_p(\lambda_n)-a
=\log_p\left(\frac{\pi(\lambda_n)}c\right)
\ne0.
\end{aligned}
\]
Thus Theorem~\ref{thm:padic-homomorphism-obstruction}  applies to this
 homomorphism $\pi_p$.
\end{proof}

\subsection{Consequences of the $p$-adic obstructions}

The additive obstruction in $\Q_p$ provides a direct proof of the following result of
Rosenblatt~\mbox{\cite[Theorem~23, p.~240]{Rosenblatt1991}}: $(2^{2^n})_{n\ge1}$ is not 
an $L^\infty$-Khinchin sequence.
We have the following more general result.

\begin{corollary}\label{cor:double-exponential-padic}
Let $p\ge 2$ be a prime and  let $(a_n)$ be a strictly increasing sequence of positive integers such that $|a_n|_p\to 0$ as $n\to \infty$. Then
$(p^{a_n})_{n\ge1}$ is not an $L^\infty$-Khinchin sequence.
\end{corollary}

\begin{proof}
The valuation $\pi(m)=v_p(m)$ is a homomorphism from $(\N,\cdot)$ to
$(\Q_p,+)$. Since $\pi(p^{a_n})=a_n\ne0$ and $a_n\to0$ in $\Q_p$,
Theorem~\ref{thm:padic-homomorphism-obstruction} applies.
\end{proof}

\begin{example}\label{exam:semigroup-exponential-padic}
Let $S(p_1, p_2, \cdots, p_d)$ be the semigroup generated by $d$ primes $p_1, p_2, \cdots, p_d$.
Let $(a_n)_{n\ge1}$ be a sequence of distinct elements of
$S(p_1,\ldots,p_d)$. For every integer $q\ge2$, the sequence
$(q^{a_n})_{n\ge1}$ is not an $L^\infty$-Khinchin sequence.
\end{example}

\begin{proof} For any $n\ge 1$, look at the vector of valuations $(v_{p_1}(a_n), v_{p_2}(a_n), \cdots, v_{p_d}(a_n))$.
Since $a_n$ tends to infinity, there exist $j\in\{1,\ldots,d\}$ and a set
$I\subseteq\N$ of upper density at least $1/d$ such that $v_{p_j}(a_n)\to\infty$ as $n\in I$ 
tends to infinity, so that
\begin{equation}\label{eq:S(p)1}
 |a_n|_{p_j}\longrightarrow 0 
 \qquad(n\to\infty,\ n\in I).
\end{equation}

Choose a prime divisor $r$ of $q$ and consider the homomorphism
\[
 \pi:(\N,\cdot)\longrightarrow(\Q_{p_j},+),
 \qquad \pi(m):=\frac{v_r(m)}{v_r(q)}.
\]
This homomorphism satisfies 
\begin{equation}\label{eq:S(p)2}
\pi(q^{a_n}) = \frac{v_r (q^{a_n})}{v_r(q)} = a_n \not=0.
\end{equation}
Thus \eqref{eq:S(p)1} and \eqref{eq:S(p)2} allow us to conclude by 
Theorem~\ref{thm:padic-homomorphism-obstruction}.
\end{proof}

\begin{corollary}\label{cor:multiplicative-p-obs}
Let $(a_n)_{n\ge1}$ be a sequence of distinct elements of
$S(p_1,\ldots,p_d)$, and let $F:\N\to\N$ be eventually strictly
increasing. Suppose that for each $1\le i\le d$ there is a
$c_i\in\Q_{p_i}^\times$ such that
\[
 F(m)\longrightarrow c_i\quad\text{in }\Q_{p_i}
 \qquad\text{as }m\to0\text{ in }\Q_{p_i}\text{ through }\N.
\]
Then, for every sufficiently large $n_0$, the sequence $(F(a_n))_{n\ge n_0}$
is not $L^\infty$-Khinchin.
\end{corollary}

\begin{proof}
As we have seen, for the sequence $(a_n)$ in $S(p_1, p_2, \ldots, p_d)$, there exist
an index $j\in\{1,\ldots,d\}$ and a set $I\subseteq\N$ of positive
upper density such that, with $p=p_j$,
\[
 a_n\longrightarrow0\quad\text{in }\Q_p
 \qquad(n\to\infty,\ n\in I).
\]
By the hypothesis on $F$, 
\[
 F(a_n)\longrightarrow c:= c_j \neq 0  \in \Q_p
 \qquad(n\to\infty,\ n\in I).
\]
Since the $a_n$ are distinct positive integers, $a_n\to+\infty$
in $\mathbb R$. The eventual strict increase of $F$ implies that
the values $F(a_n)$ are eventually pairwise distinct. In particular,
$F(a_n)=c $ for only finitely many $n$. 
 Theorem~\ref{multiplicative-p-obs} applies to
this tail with the homomorphism $m\mapsto m$ from $(\N,\cdot)$ into
$(\Q_p^\times,\cdot)$.
\end{proof}

In particular, let $P$ be a nonconstant polynomial with
$P(\N)\subseteq\N$ and $P(0)\ne0$. Interpolation gives $P\in\Q[t]$.
Thus $P(m)\to P(0)$ as $m\to0$ in each $\Q_{p_i}$. Since $P$
takes positive values on $\N$, it has positive leading coefficient
and is eventually strictly increasing. Given distinct elements of
$P(S(p_1,\ldots,p_d))$, 
Corollary~\ref{cor:multiplicative-p-obs} applies with $F=P$.
The same conclusion holds for arbitrary integer generators
$q_1,\ldots,q_d\ge2$, since
their semigroup is contained in the semigroup generated by the prime
divisors of $q_1\cdots q_d$.

\begin{example}\label{shift-Furstenberg}
For $S=S(p_1,\ldots,p_d)$, every sequence of distinct elements
of $S+1$ fails to be $L^\infty$-Khinchin. In particular, if $(s_n)_{n\ge1}$ is
the increasing enumeration of $S(2,3)$, then $(s_n+1)_{n\ge1}$
is not an $L^\infty$-Khinchin sequence.
\end{example}

\begin{proof}
Apply Corollary~\ref{cor:multiplicative-p-obs} with $F(m)=m+1$.
As $m\to0$ in each $\Q_{p_i}$, we have $F(m)\to1$.
\end{proof}

\color{black}

\section{Bernoulli multiplicative sequences are not $L^\infty$-Khinchin}
\label{sec:multiplicative-symbolic}

We shall prove Theorem~\ref{thm:random-sequence} by combining the transference
principle with a random-walk construction from Strassen's functional
law of the iterated logarithm. The construction follows the mechanism
used by Lacey--Petersen--Wierdl--Rudolph
\cite[Theorem~4 and the remarks following it]{Lacey1994} for centered
random walks.

\label{thm:strassen-flil}

Let $(\xi_n)_{n\ge1}$ be a sequence of independent identically distributed
real-valued random variables such that
$
    \mathbb E \xi_1=0,
$ and $
    \mathbb E \xi_1^2=\sigma^2
$ ($0<\sigma<\infty$).
Set
\begin{align*}
    s_n=\sum_{j=1}^n \xi_j,
    \qquad s_0=0.
\end{align*}
For $n\ge3$, let $W^\omega_n\in C([0,1])$ be the polygonal interpolation
determined by
\begin{align*}
    W^\omega_n\left(\frac{k}{n}\right)
    =
    \frac{s_k}{\sqrt{2n\log\log n}},
    \qquad 0\le k\le n.
\end{align*}
Equivalently, 
\begin{align*}
    W^\omega_n(t)
    =
    \frac{
        s_{\lfloor nt\rfloor}
        +(nt-\lfloor nt\rfloor)
        \xi_{\lfloor nt\rfloor+1}
    }
    {\sqrt{2n\log\log n}} \quad (0\le t\le1).
\end{align*}
Define
\begin{align*}
    \mathcal K_\sigma
    :=
    \left\{
        h\in C([0,1]) :
        h(0)=0,\;
        h \text{ is absolutely continuous},\;
        \int_0^1 |h'(t)|^2\,dt\le\sigma^2
    \right\}.
\end{align*}
Strassen's functional law of the iterated logarithm
(\!\!\!\cite[Theorem~3]{Strassen1964}) states that with probability one, the sequence $(W^\omega_n(\cdot))_{n\ge3}$ is relatively
compact in $C([0,1])$, endowed with the uniform norm, and its set of limit
points is exactly $\mathcal K_\sigma$.

\subsection{Sweeping out property of random walks}
We will work on $\T_b=\R/b\Z$, equipped with normalized Haar measure $\m_b$. Write
\begin{align*}
\|x\|_b:=\inf_{\ell\in\Z}|x-\ell b|
=b\|x/b\|_\T,\qquad x\in\R.
\end{align*}
We consider the random walk on $\mathbb{T}_b$ defined by $s_k \!\pmod b$.

\begin{lemma}\label{prop:Strassen-sweeping-out}
There is an event $\Omega_0$ of probability one such that, for every
$\omega\in\Omega_0$ and every $0<\varepsilon<1$, there exist $b>0$
and a measurable set $E\subset\T_b$ satisfying
\begin{equation}\label{eq:E-small}
\m_b(E)<\varepsilon
\end{equation}
and
\begin{equation}\label{eq:E-sweeps-all}
\sup_{n\ge1}\frac1n\sum_{k=1}^n
\one_E(x+s_k(\omega)\bmod b)>1-\varepsilon
\qquad\text{for every }x\in\T_b.
\end{equation}
\end{lemma}

\begin{proof}
For $n \ge3$, put
\(
L_n:=\sqrt{2n\log\log n}.
\)
Let $\Omega_0$ be the event such that $\{W_n^\omega(\cdot)\}_{n\ge 3}$
is relatively compact and has $\mathcal{K}_\sigma$ as its set of limit
points. By Strassen's theorem, $\Omega_0$ has probability $1$.

Fix $\omega\in\Omega_0$ and $0<\varepsilon<1$. Choose points
$y_1,\ldots,y_K\in[0,1)$ such that
\begin{equation}\label{eq:torus-cover}
\mathbb T
=
\bigcup_{j=1}^K
\left\{z\in\mathbb T:
\|z+y_j\|_{\mathbb T}<\frac{\varepsilon}{8}
\right\}.
\end{equation}
Then choose ordered points $0<v_1<\cdots<v_K<1$ in the interval $(0,1)$, for example $v_j = (\varepsilon/16)^{K-j}v_K$ ($0<v_K<1$ being arbitrary),   so that
\[
v_{j-1}<\frac{\varepsilon}{8}v_j
\qquad (2\le j\le K),
\]
and set $u_j=(\varepsilon/8)v_j$. Then
$[u_j,v_j]$ are pairwise disjoint intervals in $(0,1)$, and
\begin{equation}\label{eq:interval-density}
\frac{v_j-u_j}{v_j}
=1-\frac{\varepsilon}{8}
\qquad (1\le j\le K).
\end{equation}

Let $H\in C[0,1]$ be piecewise linear, with $H(0)=0$ and
\[
H(t)=y_j\qquad (u_j\le t\le v_j,\ 1\le j\le K).
\]
Extend it, for example, constantly on $[v_K,1]$. Since
$\int_0^1|H'(t)|^2\,dt<\infty$, we may choose ${c_0}>0$ small enough
that $h:={c_0} H$ belongs to $\mathcal K_\sigma$. By the choice of
$\omega\in \Omega_0$, there exist arbitrarily large integers $N$ such that
\begin{equation}\label{eq:strassen-approx}
\sup_{0\le t\le1}|W_N^\omega(t)-h(t)|
<\frac{\varepsilon{c_0}}{8}.
\end{equation}
Fix such an $N$, sufficiently large that
\[
\frac{2}{Nv_j}<\frac{\varepsilon}{8}
\qquad (1\le j\le K),
\]
and put $b:={c_0} L_N$. Define
\[
E:=
\left\{z\in\T_b:\|z\|_b\le\frac{\varepsilon b}{4}\right\}.
\]
Since $0<\varepsilon<1$,
\(
\mathfrak m_b(E)=\frac{\varepsilon}{2}<\varepsilon.
\) So, \eqref{eq:E-small} is satisfied.

If $k/N\in[u_j,v_j]$, then $h(k/N)={c_0} y_j$ and
$W_N^\omega(k/N)=s_k(\omega)/L_N$. Thus
\begin{equation}\label{eq:rk-close-cj}
|s_k(\omega)-by_j|
<\frac{\varepsilon b}{8}
\qquad\text{whenever } k/N\in[u_j,v_j].
\end{equation}

Now fix $x\in\T_b$. By \eqref{eq:torus-cover}, there is an index
$j=j(x)$ such that
\[
\|x+by_j\|_b<\frac{\varepsilon b}{8}.
\]
Together with \eqref{eq:rk-close-cj}, this gives
\[
\|x+s_k(\omega)\|_b
\le \|x+by_j\|_b+|s_k(\omega)-by_j|
<\frac{\varepsilon b}{4}
\]
whenever $k/N\in[u_j,v_j]$. Hence
\begin{equation}\label{eq:block-in-E}
x+s_k(\omega)\pmod b\in E
\qquad\text{whenever } k/N\in[u_j,v_j].
\end{equation}

Set $n_j:=\lfloor Nv_j\rfloor$. At least
$N(v_j-u_j)-2$ indices $1\le k\le n_j$ satisfy
$k/N\in[u_j,v_j]$. Since $n_j\le Nv_j$, 
\eqref{eq:block-in-E} and \eqref{eq:interval-density} yield
\begin{align*}
\frac1{n_j}\sum_{k=1}^{n_j}
\one_E(x+s_k(\omega)\bmod b)
&\ge \frac{N(v_j-u_j)-2}{Nv_j}\\
&=1-\frac{\varepsilon}{8}-\frac{2}{Nv_j}\\
&>1-\frac{\varepsilon}{4}
>1-\varepsilon.
\end{align*}
The index $j$ may depend on $x$, but the same $b$ and $E$ work for
every $x\in\T_b$. This proves \eqref{eq:E-sweeps-all}.
\end{proof}

The event $\Omega_0$ is chosen once from Strassen's theorem, before
$\varepsilon$ is chosen. The period $b$ and the set $E$ may depend on
both $\omega$ and $\varepsilon$. 

For the Bernoulli application, we choose centered increments taking
only two values. Pairwise coprimality of the generators allows us
to realize these increments by a homomorphism from $\N$ to $\R$.

\subsection{Proof of Theorem~\ref{thm:random-sequence}: Bernoulli random multiplicative sequences}

{We shall need the following simple lemma.}
\begin{lemma}\label{lem:homom}
Let $(q_n)_{n\ge1}$ be pairwise coprime integers with $q_n\ge2$, and let
$(a_n)_{n\ge1}\subset\R$ be arbitrary. Then there exists a semigroup
homomorphism
$
\pi:(\N,\cdot)\longrightarrow (\R,+)
$
such that
\begin{align*}
\pi(q_n)=a_n
\qquad\text{for every }n\ge1.
\end{align*}
\end{lemma}

\begin{proof}
For each $n$, choose a prime $p_n\mid q_n$. Since the $q_n$ are pairwise
coprime, the primes $p_n$ are distinct. Define  $(c_p)$ on the set $\mathcal{P}$ of 
primes by
\begin{align*}
c_{p_n}
=
\frac{a_n}{v_{p_n}(q_n)} \ \text{for}  \ n\ge 1, \qquad c_p=0 \ \text{for} \  p \in \mathcal{P}\setminus\{p_n: n \ge 1\}.
\end{align*}
For $m\in\N$, define
\begin{align*}
\pi(m)
=
\sum_{p} v_p(m)c_p.
\end{align*}
The sum is finite. Since
$
v_p(mk)=v_p(m)+v_p(k)
$,  $\pi$ is a semigroup homomorphism. 
Moreover, by pairwise coprimality, $v_{p_m}(q_n) \not =0$ if and only if $m= n$,  hence
\begin{align*}
\pi(q_n)
=
v_{p_n}(q_n)c_{p_n}
=
a_n.
\end{align*}
\end{proof}

\begin{proof}[Proof of Theorem~\ref{thm:random-sequence}]
Since the support of $\rho$ has at least two points, choose $r\in\operatorname{supp}\rho$ so that
$\theta:=\rho(\{r\})\in(0,1)$.  We construct a centered random walk on $\mathbb{R}$  with two possible increments. 
Put 
\begin{equation}\label{eq:centered-increments}
a_j:=\begin{cases}
\theta-1,&j=r,\\
\theta,&j\in \mathbb{N}\setminus \{r\}.
\end{cases}
\end{equation}
Then $\xi_k(\omega):=a_{\omega_k}$ are i.i.d. under $\rho^\N$, with
\begin{align*}
\Ebb\xi_1&=\theta(\theta-1)+(1-\theta)\theta=0,\\
\operatorname{Var}(\xi_1)
&=\theta(\theta-1)^2+(1-\theta)\theta^2
=\theta(1-\theta)>0.
\end{align*}

 Now choose a homomorphism
$\ell:(\N,\cdot)\to(\R,+)$ satisfying 
$$
\ell(q_n)=a_n, \quad \text{for} \ n \ge 1.
$$
Since the $q_n$ are pairwise  coprime, such a homomorphism exists by Lemma~\ref{lem:homom}.
 Then 
\begin{align*}
\ell(\lambda_k)=\ell(q_{\omega_1}\cdots q_{\omega_k})
=s_k(\omega).
\end{align*}

Let $\Omega_0$ be the probability-one event appearing in Lemma~\ref{prop:Strassen-sweeping-out}, and 
fix $\omega\in\Omega_0$.  Choose
$0<\varepsilon_m<1$ with $\varepsilon_m\to0$.
For each $m$, Lemma~\ref{prop:Strassen-sweeping-out} provides 
a positive number $b_m$ and a Borel set
$E_m\subset\T_{b_m}$ such that
\begin{equation}\label{eq:E-small-m}
\m_{b_m}(E_m)<\varepsilon_m
\end{equation}
and
\begin{equation}\label{eq:E-sweeps-all-m}
\sup_{n\ge1}\frac1n\sum_{k=1}^n
\one_{E_m}(x+s_k(\omega)\bmod b_m)>1-\varepsilon_m
\qquad(\forall x\in\T_{b_m}).
\end{equation}
For each $m$,  we consider the  measure-preserving {$(\N,\cdot)$-action} on $\T_{b_m}$ 
defined by translations induced by the homomorphism $\ell$:
\begin{align*}
{\gamma}_m^n x:=x+\ell(n)\pmod{b_m},\qquad n\in\N.
\end{align*}
For the sequence $\Lambda=\Lambda(\omega)$, the inequalities
\eqref{eq:E-small-m} and \eqref{eq:E-sweeps-all-m} give
\begin{equation}\label{eq:random}
\mathcal D_{\varepsilon_m}(\Lambda^{{\gamma}_m,*})
\ge\|\Lambda^{{\gamma}_m,*}\one_{E_m}\|_1
\ge1-\varepsilon_m.
\end{equation}
By \eqref{eq:E-small-m} and \eqref{eq:random},
the uniform condition in
Theorem~\ref{thm:homomorphism-criterion} (ii) is violated. Hence
$\Lambda(\omega)$ is not $L^\infty$-Khinchin for every
$\omega\in\Omega_0$. 
The fact that the circle models depend on the fixed realization $\omega$ causes no difficulty, since the supremum in Theorem~\ref{thm:homomorphism-criterion} (ii) ranges over all actions.
\end{proof}

By \cite[Theorem~4.5]{FFQQ},
for every fixed $f\in L^1(\T)$,
\begin{equation}\label{eq:asae}
\text{a.s.} \ \   \m \text{-a.e.} \quad  \lim_{N\to\infty} \frac{1}{N}\sum_{n=1}^N f(q_{\omega_1}q_{\omega_2}\cdots q_{\omega_n} x) = \int_\T fd \m.
\end{equation}
This follows from the ergodicity of  the associated measure-preserving skew product. See also \cite{UlamVonNeumann1945}  for
the random ergodic theorem. 
Theorem~\ref{thm:random-sequence} shows, however, that the almost-sure event involved in \eqref{eq:asae} cannot  in general be chosen independently of $f$.

\bibliographystyle{plain}

\bibliography{KP2}

\end{document}